\documentclass[lettersize,journal]{IEEEtran}

\ifCLASSOPTIONcompsoc
\usepackage[nocompress]{cite}
\else
\usepackage{cite}
\fi

\usepackage{enumerate}
\usepackage{latexsym}
\usepackage{amsfonts}
\usepackage{amsmath}
\usepackage{amssymb}
\usepackage{color}
\usepackage{colortbl}
\usepackage{epsfig}
\usepackage{xspace}
\usepackage{graphicx}
\usepackage{arcs}
\usepackage{subfigure}
\usepackage{balance}
\usepackage[english]{babel}
\usepackage{makecell}
\usepackage{comment}
\usepackage{epsfig}
\usepackage{multirow}
\usepackage{url}
\usepackage{CJK}

\newtheorem{lemma}{Lemma}
\newtheorem{theorem}{Theorem}

\newtheorem{proposition}{Proposition}
\newtheorem{remark}{Remark}
\newtheorem{corollary}{Corollary}
\newtheorem{claim}{Claim}

\newenvironment{proof}{
	\vspace{1ex}
	{\noindent\bf Proof:}}{\eop\vspace{1ex}}
\newcommand{\eop}{\hspace*{\fill}\mbox{$\Box$}}     

\newcommand{\pll}{\kern 0.40em/\kern -0.82em /\kern 0.40em}

\begin{document}
\begin{CJK*}{GBK}{kai}

\title{Asymptotic  critical transmission  radii in random geometry graphs over  three-dimensional   regions}

\author{Jie Ding,
	Shuai Ma, 
	Xiang Wei,
	Xiaohua Xu,
	and Xinshan Zhu
	
	\IEEEcompsocitemizethanks{\IEEEcompsocthanksitem J. Ding  is with School of Computer, Jiangsu University of Science and Technology, Zhenjiang 212100, China.		\protect\\
		E-mail: jieding@just.edu.cn  

		\IEEEcompsocthanksitem  S. Ma  is with SKLSDE Lab,
		Beihang University, Beijing 100191, China.\protect\\
		E-mail: mashuai@buaa.edu.cn  
		
	 \IEEEcompsocthanksitem  X. Wei  is with School of Computer, Jiangsu University of Science and Technology, Zhenjiang 212100, China. 		 \protect\\
		E-mail: 274759197@qq.com

		\IEEEcompsocthanksitem  X. Xu  is with School of Information Engineering,
		Yangzhou University, Yangzhou 225009, China. 		\protect\\
		E-mail:  xhxu@yzu.edu.cn 
		
		\IEEEcompsocthanksitem X. Zhu is with   School of Electrical and Information Engineering, Tianjin University, Tianjin 300072, China.\protect\\
		E-mail: xszhu@tju.edu.cn
		
	}
}

{}

\maketitle

\begin{abstract}
	
	This article presents  the precise asymptotical distribution of two types of critical transmission radii, defined in terms of
	$k-$connectivity and the minimum vertex degree, for   random geometry graphs distributed over   three-dimensional   regions.

\end{abstract}


\section{Introduction and main results}\label{sec-intro}

 Three-dimensional (3D) self-organising communication networks, ranging from transitional Internet of Things (IoTs) systems including  wireless  sensor networks or Ad-Hoc networks,	to  emerging machine-type communication networks designed for unmanned aerial vehicles (UAVs), unmanned ground vehicles (UGVs), and unmanned underwater vehicles  (UUVs),  have seen widespread adoption and rapid development.
  Random geometric graphs  are a basic model of such  self-organising communication networks~\cite{Ding2025-IEEETON, Structure-1DRGG-IT2023, PJWan-IT-asymptotic-radius, Zerolaw-1DRGG-IT2009,  Haenggi-book2012, Penrose-RGG-book}, where  the devices distributed on   region $\Omega\subset\mathbb{R}^3$ are naturally modeled as the vertices of a graph, which  can be represented by a point process $\chi_n$   uniformly distributed on $\Omega$.     Each devices are  assumed to have the same transmission range (or radius) $r_n$,  so  the communication link establashes  between any pair of devices with   their Euclidean distance   less than $r_n$. These links naturally form the edges of the graph. Such obtained random geometry graph is denoted by  $G(\chi_n,r_n)$.

Connectivity, playing a crucial role in robust network  communications, is
 an interesting and  important topological property for  a    graph. If there does not exist a set of $k-1$ vertices  whose removal disconnects a graph, then this graph is said to be  $k-$connected.  Let  $\kappa$  denote the largest number  $k$ such that a graph is $k-$connected, and random variable
$\rho(\chi_n;\kappa\geq k)$ be the minimum $r_n$ such that random geometry graph
$G(\chi_n,r_n)$ is $k-$connected.  In addition,  the minimum vertex degree of a graph is denoted by $\delta$, and notation $\rho(\chi_n;\delta\geq k)$ denotes the minimum $r_n$ such that   $G(\chi_n,r_n)$ has the smallest vertex degree $k$.  As pointed out in~\cite{Ding2025-IEEETON},  the performance metrics of  network communications, such like coverage,  throughput, routing efficiency and energy consumption, rely on the critical transmission   radius $\rho(\chi_n;\kappa\geq k)$ or $\rho(\chi_n;\delta\geq k)$.

\begin{table*}
	\begin{center}\caption{Related works}\label{table:relatec-work}
		\begin{tabular}{|c|c|c|c|}
			\hline\hline
			Region $\Omega$ & $r_n$ &   Result  & Reference   \\ \hline
			Square $\subset \mathbb{R}^2$   & $r_n=\sqrt{\frac{\log n+c}{\pi n}}$
			& $ \Pr\left\{\rho(\chi_n;\delta\geq 1)\leq r_n\right\}\sim \exp\left(-e^{-c}\right)$
			& Dette and Henze~\cite{DH-largest-NN}, 1989  \\\hline
			
			Disk $\subset\mathbb{R}^2$  & $r_n=\sqrt{\frac{\log n+c}{\pi n}}$
			& $ \Pr\left\{\rho(\chi_n;\kappa\geq 1)\leq r_n\right\}\sim \exp\left(-e^{-c}\right)$
			& Gupta and Kumar~\cite{GK-Critical-Power}, 1998  \\\hline
			
			Square $\subset \mathbb{R}^d$  &  $\frac n{k!}\int_{\Omega}\left(n|B(x,r_n)\cap\Omega|\right)^k e^{-n|B(x,r_n)\cap\Omega|}\mathrm{d}x\sim \lambda$
			& \makecell{$ \Pr\left\{\rho(\chi_n;\delta\geq k+1)\leq r_n\right\}\sim \exp\left(-\lambda\right)$\\
				$ \Pr\left\{\rho(\chi_n;\kappa\geq k+1)\leq r_n\right\}\sim \exp\left(-\lambda\right)$ }
			& Penrose~\cite{Penrose-RGG-book}, 2003  \\\hline
			
			Square$\subset \mathbb{R}^2$
			& \makecell{ $r_n=\sqrt{\frac{\log n+(2k-1)\log\log n+\xi}{\pi n }}$, with  \\
				$ \xi=\left\{\begin{array}{cc}
					-2\log \left(\sqrt{e^{-c}+\frac{\pi}{4}}-\frac{\sqrt{\pi}}{2}\right), & k=1,\\
					2\log\frac{\sqrt{\pi}}{2^{k-1}k!}+2c, & k>1, \\
				\end{array}
				\right. $ }
			
			& \makecell{$ \Pr\left\{\rho(\chi_n;\delta\geq k+1)\leq r_n\right\}\sim \exp\left(-e^{-c}\right)$\\
				$ \Pr\left\{\rho(\chi_n;\kappa\geq k+1)\leq r_n\right\}\sim \exp\left(-e^{-c}\right)$ }
			& Wan, \emph{et.al.}~\cite{PJWan-IT-asymptotic-radius}, 2010
			 \\\hline

			Disk $\subset \mathbb{R}^2$
			& \makecell{ $r_n=\sqrt{\frac{\log n+(2k-1)\log\log n+\xi}{\pi n }}$, with  \\
				$ \xi=\left\{\begin{array}{cc}
					-2\log \left(\sqrt{e^{-c}+\frac{\pi^2}{16}}-\frac{\pi}{4}\right), & k=1, \\
					2\log\frac{{\pi}}{2^{k}k!}+2c, & k>1. \\
				\end{array}
				\right.$}
			
			& \makecell{$ \Pr\left\{\rho(\chi_n;\delta\geq k+1)\leq r_n\right\}\sim \exp\left(-e^{-c}\right)$\\
				$ \Pr\left\{\rho(\chi_n;\kappa\geq k+1)\leq r_n\right\}\sim \exp\left(-e^{-c}\right)$ }
			& Wan, \emph{et.al.}~\cite{ PJWan-IT-asymptotic-radius}, 2010 \\\hline

			\makecell{Convex region\\ $\subset \mathbb{R}^2$}
			& $ r_n=\sqrt{\frac{\log n+c}{\pi n }}$
			
			& \makecell{$ \Pr\left\{\rho(\chi_n;\delta\geq 1)\leq r_n\right\}\sim \exp\left(-e^{-c}\right)$\\
				$ \Pr\left\{\rho(\chi_n;\kappa\geq 1)\leq r_n\right\}\sim \exp\left(-e^{-c}\right)$ }
			& Ding, \emph{et.al.}~\cite{Ding2025-IEEETON}, 2025 \\\hline

			\makecell{Convex region\\ $\subset \mathbb{R}^2$, with \\ $|\partial \Omega|=l$}
			& \makecell{ $r_n=\sqrt{\frac{\log n+(2k-1)\log\log n+\xi}{\pi n }}$, with  \\
				$ \xi=\left\{\begin{array}{cc}
					-2\log \left(\sqrt{e^{-c}+\frac{\pi l^2}{64}}-\frac{l\sqrt{\pi}}{8}\right), & k=1, \\
					2\log \left(\frac{l\sqrt{\pi}}{2^{k+1}k!}\right)+2c, & k>1. \\
				\end{array}
				\right. $}
			
			& \makecell{$ \Pr\left\{\rho(\chi_n;\delta\geq k+1)\leq r_n\right\}\sim \exp\left(-e^{-c}\right)$\\
				$ \Pr\left\{\rho(\chi_n;\kappa\geq k+1)\leq r_n\right\}\sim \exp\left(-e^{-c}\right)$ }
			& Ding, \emph{et.al.}~\cite{Ding2025-IEEETON}, 2025 \\\hline

		\end{tabular}
	\end{center}
\end{table*}

Dette and Henze~\cite{DH-largest-NN} drived the precise asymptotic distribution of $\rho(\chi_n;\delta\geq 1)$ for a unit-area square on $\mathbb{R}^2$.   For the case of unit-area disk   on $\mathbb{R}^2$,
 Gupta and Kumar gives the precise
 asymptotic distribution of $\rho(\chi_n;\kappa\geq 1)$ in~\cite{GK-Critical-Power}.  Penrose~\cite{Penrose97:longest-edge} proved that the critical transmission radius in terms of connectivity equals   the longest edge of the minimal spanning tree, and further gave  the asymptotic distribution of the longest edge in~\cite{Penrose99:Strong-law-longest-edge}.
In addition,  Wang~\cite{XBWang2013:MotionCast-book} defined  a type of $(k,m)-$connectivity and derived  the precise asymptotic distribution of the critical transmission range over a square region.

For a hypercube $\Omega\subset\mathbb{R}^d\;(d\geq 2)$,
Penrose     demonstrated the following result:
\renewcommand*{\thetheorem}{\Alph{theorem}}
\begin{theorem}\label{thm:Penrose-formulae}(\cite{Penrose-RGG-book,Penrose-k-connectivity})
	Let $k>0,\lambda\in \mathbb{R}^+$,  $\Omega\subset\mathbb{R}^d\;(d\geq 2)$ be a unit-volume hypercube.
	If sequence $\left(r_n\right)_{n\geq 1}$ satisfies
	\begin{equation}\label{eq:formula-0}
		\frac n{k!}\int_{\Omega}\left(n|B(x,r_n)\cap\Omega|\right)^k e^{-n|B(x,r_n)\cap\Omega|}\mathrm{d}x\sim \lambda,
	\end{equation}
	then as $n\rightarrow\infty$, the probabilities of the two events
	$\rho(\chi_n;\delta\geq k+1)\leq r_n$ and $\rho(\chi_n;\kappa\geq
	k+1)\leq r_n$ both converge to $e^{-\lambda}$.
\end{theorem}

Further,    Wan \emph{et al.}~\cite{PJWan-IT-asymptotic-radius} and~\cite{PWan04:Asymptotic-critical-transmission-MobileHoc}  provided the explicit forms of $r_n$ satisfying   formula
(\ref{eq:formula-0}) in the following two theorems.
\begin{theorem}(\cite{PJWan-IT-asymptotic-radius, PWan04:Asymptotic-critical-transmission-MobileHoc})
(1). If $\Omega\subset \mathbb{R}^2$ is a unit-area square. Let
	$$
	r_n=\sqrt{\frac{\log n+(2k-1)\log\log n+\xi}{\pi n }},
	$$
	where $\xi$ satisfies
	$$ \xi=\left\{\begin{array}{cc}
		-2\log \left(\sqrt{e^{-c}+\frac{\pi}{4}}-\frac{\sqrt{\pi}}{2}\right), & k=1,\\
		2\log\frac{\sqrt{\pi}}{2^{k-1}k!}+2c, & k>1, \\
	\end{array}
	\right.
	$$
	then   formula (\ref{eq:formula-0}) holds with $\lambda=e^{-c}$.

\par (2). If  $\Omega\subset \mathbb{R}^2$ is a unit-area disk. Let
	$$
	r_n=\sqrt{\frac{\log n+(2k-1)\log\log n+\xi}{\pi n }},
	$$
	where $\xi$ satisfies
	$$ \xi=\left\{\begin{array}{cc}
		-2\log \left(\sqrt{e^{-c}+\frac{\pi^2}{16}}-\frac{\pi}{4}\right), & k=1, \\
		2\log\frac{{\pi}}{2^{k}k!}+2c, & k>1. \\
	\end{array}
	\right.
	$$
	then   formula~(\ref{eq:formula-0}) holds  with $\lambda=e^{-c}$.

In the above two cases,  the probabilities of the two events
$\rho(\chi_n;\delta\geq k+1)\leq r_n$ and $\rho(\chi_n;\kappa\geq
k+1)\leq r_n$ both converge to $\exp(-e^{-c})$ as $n$ tends to
infinity.
	
\end{theorem}

Then, these results are  generalised into the following
\renewcommand*{\thetheorem}{\Alph{theorem}}
\begin{theorem}\label{thm:combined-Main}(\cite{Ding2025-IEEETON})
	Let $\Omega\subset \mathbb{R}^2$ be a unit-area convex region with  the length of
	the boundary $\mathrm{Len}(\partial\Omega)$,  $k\geq 0$ be an integer and $c>0$ be a constant.
	\par (i) If $k>0$, let
	\begin{equation*}\label{eq:Theorem-formula-1}
		r_n=\sqrt{\frac{\log n+(2k-1)\log\log n+\xi}{\pi n }},
	\end{equation*}
	where $\xi$ satisfies
	$$ \left\{\begin{array}{cc}
		\xi=-2\log \left(\sqrt{e^{-c}+\frac{\pi (\mathrm{Len}(\partial\Omega))^2}{64}}-\frac{\mathrm{Len}(\partial\Omega)\sqrt{\pi}}{8}\right), & k=1, \\
		\xi=2\log \left(\frac{\mathrm{Len}(\partial\Omega)\sqrt{\pi}}{2^{k+1}k!}\right)+2c, & k>1. \\
	\end{array}
	\right.
	$$
	\par (ii) If $k=0$, let $r_n=\sqrt{\frac{\log n+c}{\pi n }}.$

   Then  formula~(\ref{eq:formula-0}) holds  with $\lambda=e^{-c}$,
	and therefore,  the probabilities of the two events
	$\rho(\chi_n;\delta\geq k+1)\leq r_n$ and $\rho(\chi_n;\kappa\geq
	k+1)\leq r_n$ both converge to $\exp\left(-e^{-c}\right)$  as $n\rightarrow\infty$.
\end{theorem}
Theorem~\ref{thm:combined-Main}  clearly reveals how the region shape impacts on the critical transmission ranges.  The above mentioned works are listed in Table~\ref{table:relatec-work} for comparison.

The main contribution of this article is to establish explicit expressions for the critical radii---defined in terms of $k$-connectivity and minimum vertex degree---and to demonstrate their asymptotic distributions, for the random geometric graphs distribued  over general three-dimensional regions in $\mathbb{R}^3$.
Moreover, although our proof of the main results follows the framework developed in~\cite{Ding2025-IEEETON} for two-dimensional domains, this article introduces more delicate boundary-treatment techniques to establish the corresponding results in three dimensions.

 For convenience, throughout this article, we define
\begin{equation}\label{eq:psi-function}
	\psi^k_{n,r}(x)= \frac{\left(n|B(x,r)\cap\Omega|\right)^k
		e^{-n|B(x,r)\cap\Omega|}}{k!}.
\end{equation}
Function $\psi^k_{n,r}(x)$  indicates the probability of that the node at $x$ has $k$ degree, for a Poisson point process  which will be introduced  in Section~\ref{sec:Pro-2}.
Our main result is
\renewcommand*{\thetheorem}{\arabic{theorem}}
 \setcounter{theorem}{0}
\begin{theorem}\label{thm:Main3D}
	Suppose $\Omega\subset \mathbb{R}^3$ is a simply connected compact region with unit-volume. Point process $\chi_n$ is  uniformly distributed on $\Omega$.
The boundary $\partial\Omega$ of the region has finite area $\mathrm{Area}(\partial\Omega)$ and
satisfies  one of the following conditions:
\begin{description}
\item[$C_{\mathrm{I}}$]: The boundary $\partial\Omega$ is $C^2$ smooth everywhere.
\item[$C_{\mathrm{II}}$]: The boundary $\partial\Omega$ is piecewise $C^2$ smooth, and the non-smooth points lie on  a curve $\sigma$ of finite   length and zero area.
The maximum  and minimum principal curvatures  of surface  $\partial\Omega\backslash \sigma$  are uniformly bounded.
In addition, $\liminf_{r\rightarrow 0} \inf_{x\in\Omega}\frac{\mathrm{Vol}(B(x,r)\cap \Omega )}{\mathrm{Vol}(B(x,r))}>\frac{1}{4}.$
\item[$C_{\mathrm{III}}$]: The boundary $\partial\Omega$ is $C^2$ smooth everywhere, except for a finite number of points $\{A_i\}_{i=1}^M$.
The maximum  and minimum principal curvatures  of surface  $\partial\Omega\backslash \{A_i\}_{i=1}^M$  are uniformly bounded.
 Additionally, the region satisfies $\liminf_{r\rightarrow 0} \inf_{x\in\Omega}\frac{\mathrm{Vol}(B(x,r)\cap \Omega )}{\mathrm{Vol}(B(x,r))}>0.$
    \end{description}
 Assume $k\geq 0$ is an integer and $c>0$ is a constant.
      Let
	\begin{equation}\label{eq:Theorem-radius}
		r_n= \left(\frac{\log n+(\frac{3k}{2}-1)\log\log n+\xi}{  \pi n  }\right)^{\frac13},
	\end{equation}
	where   $\xi$  solves
	\begin{equation}\label{eq:Theorem-radius-1}
		\mathrm{Area}(\partial\Omega)\frac{e^{-\frac{2\xi}{3}}}{ \pi^{\frac{1}{3}}} \left(\frac{2}{3}\right)^k\frac{1}{k!} =e^{-c},
	\end{equation}
	or equivalently,
	\begin{equation}\label{eq:Theorem-radius-2}
		c=-\log(\mathrm{Area}(\partial\Omega))-k\log\left(\frac{2}{3}\right)+\log(k!) +\frac{2\xi}{3}+\frac{1}{3}\log\left(\pi\right).
	\end{equation}
	Then  	$		
	\lim_{n\rightarrow\infty}  \int_{\Omega}n\psi^k_{n,r}(x)\mathrm{d}x= e^{-c}, $
	and therefore, the probabilities of the two events
	$\rho(\chi_n;\delta\geq k+1)\leq r_n$ and $\rho(\chi_n;\kappa\geq
	k+1)\leq r_n$ both converge to $\exp\left(-e^{-c}\right)$  as $n\rightarrow\infty$.

\end{theorem}

 Thoughout this article, a surface is $C^2$ smooth means that the function of  the surface has second-order continuous derivatives everywhere.  If surface $\partial\Omega$ is not $C^2$ smooth  everywhere, as
 Condition~$C_{\mathrm{II}}$ states,  the singularity part of the boundary  (e.g., the edges and vertices of a polyhedron) is contained in $\sigma$ and  $\partial\Omega\backslash \sigma$ is   $C^2$ smooth.  For convenience and without the loss of generality,  here we may assume that      $\sigma\subset\partial\Omega$ is a curve of finite length.
For example,  the edges and vertices of a polyhedron can    be connected by a single curve, i.e., they can lie on a single curve with zero area and finite length.

It is clearly to see that there is a fundamental difference between the three-dimensional and two-dimensional cases.
In the two-dimensional scenario (see Theorem~\ref{thm:combined-Main}), the length of the region boundary has no impacts whatsoever on  $1-$connectivity and  $1-$minimal vertex degree. In contrast, for the three-dimensional case, as Theorem~\ref{thm:Main3D} shows,   the boundary does influence both $1-$connectivity and  $1-$minimal vertex degree through formulas (\ref{eq:Theorem-radius})-(\ref{eq:Theorem-radius-2}).

Hereafter, we always assume that  $\Omega  \subset \mathbb{R}^3$ is  simply connected, with boundary $\partial\Omega$ of area $\mathrm{Area}(\partial\Omega)$.
The following remark gives some estimation about the sequence $r_n$, which is straightforward and will be often used.
\begin{remark}\label{rem:remark3}
	Sequence $r_n$   given by (\ref{eq:Theorem-radius})-(\ref{eq:Theorem-radius-2}) in Theorem~\ref{thm:Main3D} satisfies
	that \\
	(i).   $ {n \pi r^3}= \left[\log n+\left(\frac{3k}{2}-1\right)\log\log n +\xi\right] $ 
and
  $$
	\frac{n}{k!}\left( \frac{4 }{3}  n\pi r^3\right)^{k}e^{-\frac{4 }{3}n\pi r^3}
	\sim  \frac{n(\frac{4}{3}\log n+o(\log n))^{k}}{k!e^{(\frac{4}{3}\log n+o(\log n))}}=o(1).
	$$
	
   (ii).  If constant $C^*>0$, then	
	$$
	\left(   n\pi r^3\right)^{k+1}e^{-\frac{4 C^*}{3}n\pi r^3}
	\sim   \frac{( \log n+o(\log n))^{k+1}}{e^{ \frac{4C^*}{3}\log n+o(\log n)}}=o(1).
	$$
	
	(iii).  If constant $C_*>\frac{1}{4}$, then
	$$
    \frac{ \left(   n\pi r^3\right)^{k+1}}{r}e^{-\frac{4 C_*}{3}n\pi r^3}
	=  \frac{( \log n+o(\log n))^{k+\frac{2}{3}}} { \Theta(1) n^{\frac{4C_*-1}{3}}  (\log n)^{\frac{4C_*}{3}\frac{3k-2}{2}}}=o(1).
	$$

\end{remark}

The proof  of Theorem~\ref{thm:Main3D} follows the framework presented in~\cite{Ding2025-IEEETON}: to prove Theorem~\ref{thm:Main3D} it suffices to prove the following three propositions.

\begin{proposition}\label{pro:Explicit-Form}
Under the assumptions of Theorem~\ref{thm:Main3D},
	\begin{equation}\label{eq:Explicit-Form}
			\lim_{n\rightarrow\infty} \int_{\Omega}n\psi^k_{n,r}(x)\mathrm{d}x= e^{-c}.
	\end{equation}
\end{proposition}

\begin{proposition}\label{pro:conclusion-1}   Under the assumptions of Theorem~\ref{thm:Main3D},
	\begin{equation}\label{eq:conclusion-1}
		\lim_{n\rightarrow\infty}\Pr\left\{\rho(\chi_n;\delta\geq k+1)\leq r_n\right\} =\exp(-e^{-c}).
	\end{equation}
\end{proposition}

\begin{proposition}\label{pro:conclusion-2}  Under the assumptions of Theorem~\ref{thm:Main3D},
	\begin{equation}\label{eq:conclusion-2}
		\lim_{n\rightarrow\infty}\Pr\left\{\rho(\chi_n;\delta\geq k+1)=\rho(\chi_n;\kappa\geq
		k+1)\right\}=1.
	\end{equation}
\end{proposition}

  We should point out that  the  conclusion also holds when the uniform point process $\chi_n$
  in Theorem~\ref{thm:Main3D} is replaced by a  Poisson point process $\mathcal{P}_n$.
  In fact, when proving Theorem~\ref{thm:Main3D}, we first demonstrate that its conclusion holds for a Poisson point process, and then use the de-Poissonization technique to show that it also holds for the  $\chi_n$   point process (see Section~\ref{sec:Pro-2} for details).

In the following, Proposition~\ref{pro:Explicit-Form}-\ref{pro:conclusion-2} will be proved in Section~\ref{sec:Pro-1} to Section~\ref{sec:Pro-3} respectively, while some preliminary results on regions in $\mathbb{R}^3$   are presented in Section~\ref{sec:Preliminary-Results} which will be used for the proofs.

\par
We use the following notations throughout this article.
(1)   $B(x,r)\subset \mathbb{R}^3$ is a ball or sphere
centered at $x$ with  radius $r$.
(2) Notation $|A|$ is a short  for the volume of a measurable set $A\subset \mathbb{R}^3$ and $\|\cdot\|$ represents the length of a line segment. $\mathrm{Area}(\cdot)$ denotes the area of a surface.
(3) $\mathrm{dist}(x, \Sigma)=\inf_{y\in \Sigma}\|xy\|$ where $x$ is a point and $\Sigma$ is a set.
(4) For any two nonnegative functions $f(n)$ and $g(n)$, if there exist two constants $0<c_1<c_2$ such that
$c_1g(n)\leq f(n)\leq c_2g(n)$ for any sufficiently large $n$, then  denote $f(n)=\Theta(g(n))$.
We also use notations $f(n)=o(g(n))$ and $f(n)\sim g(n)$ to denote that $\lim\limits_{n\rightarrow\infty}\frac{f(n)}{g(n)}=0$ and $\lim\limits_{n\rightarrow\infty}\frac{f(n)}{g(n)}=1$, respectively.
The notations introduced here and to appear later are  summarised in a  Table~\ref{table:notation}.

\begin{table*}
	\begin{center}\caption{Notations and descriptions}\label{table:notation}
		\begin{tabular}{|c|c|}
			\hline\hline
			Notations & Descriptions \\ \hline
			$\Omega$ & A compact and simply connected    region in $\mathbb{R}^3$ with unit-volume \\\hline
		
			$\partial\Omega$ &  The boundary of region $\Omega$ \\\hline
		
			$ \mathrm{dist}(x, \Sigma)$ 			& $ \mathrm{dist}(x,\Sigma)=\inf_{y\in \Sigma}\|xy\|$, the distance of point $x$ to set $\Sigma$. \\\hline

			$W_x$ &  The nearest point on $\partial\Omega$ to point $x$, with $\|xW_x\|=\mathrm{dist}(x,\partial\Omega)$ \\\hline
			
			$\mathcal{T}_{W_x}(\partial\Omega),\mathcal{T}_{W_x}$ &  The tangent plane of $\partial\Omega$ at point $W_x$ \\\hline
			$\pll$ &  parallel between planes \\\hline	
			
			$\gamma_{P}(\mathcal{T}_{W_x})$ &  The plane which passes through point $P$ and is parallel to $\mathcal{T}_{W_x}(\partial\Omega)$\\\hline

			$B(x,r)$   & Sphere (ball) centered at $x$ with radius $r$ \\\hline					
			$B(x,r)|_{\gamma_{E}(\mathcal{T}_{W_x})}$ & The major segment of sphere $B(x,r)$ cut by plane $\gamma_{E}(\mathcal{T}_{W_x})$ \\\hline


			$\mathrm{C_{max}}(\partial\Omega|_O)$,	$\mathrm{C_{min}}(\partial\Omega|_O)$ &  Maximum  and minimum principal curvatures  of surface  $\partial{\Omega}$ at point $O$ respectively\\\hline

            $G(\partial\Omega), G(\partial\Omega\backslash\sigma)$ &  Curvature  functions of surfaces $\partial\Omega$ and $\Omega\backslash\sigma$ respectively (see Lemma~\ref{lemma:G(Omega)} and Corollary~\ref{corollary:Bounds-B(x,r)-Omega-sigma})  \\\hline

			$\chi_n$ &   A uniform $n$-point process over $\Omega$   \\\hline
			$ \mathcal{P}_n$ &  A homogeneous Poisson process of intensity $n$ (i.e., $n|\Omega|$) on $\Omega$  \\\hline
			$G(\chi_n,r_n)$ & A graph obtained by connecting each pair of points in $\chi_n$ if their distance
			is less than $r_n$   \\\hline
			
			$G(\mathcal{P}_n,r_n)$ & A graph obtained by connecting each pair of points in $\mathcal{P}_n$ if their distance
			is less than $r_n$    \\\hline
			
			$\rho(\chi_n;\kappa\geq k)$  &   The minimum $r_n$ such that
			$G(\chi_n,r_n)$ is $k-$connected   \\\hline
			$\rho(\chi_n;\delta\geq k)$ &  The minimum $r_n$ such that $G(\chi_n,r_n)$ has minimum degree $k$         \\\hline
			$\rho(\mathcal{P}_n;\delta\geq k)$& The minimum $r_n$ such that $G( \mathcal{P}_n,r_n)$ has minimum degree $k$ \\\hline

			$|\cdot|$ & Volume of a measurable set in $\mathbb{R}^3$  \\\hline
			$\mathrm{Area}(\cdot)$ & Area of a measurable set in $\mathbb{R}^3$  \\\hline
			$\|\cdot\|$   & Length of a line segment \\\hline

			$ \psi^k_{n,r}(x)$ & Short for $\frac{1}{k!}(n|B(x,r)\cap \Omega|)^k\exp(-n|B(x,r)\cap \Omega|)$  \\\hline
			
			$t(x,r), t(x)$  &  $	t(x,r)=\inf\{\mathrm{dist}(x,\gamma_{E'}(\mathcal{T}_{W_x}))\mid E'\in\partial\Omega\cap\partial B(x,r)   \}$, see Lemma~\ref{lemma:Bounds-B(x,r)-Omega} and (\ref{eq:t(x,r)})\\\hline

			$a(r,t)$, $a(t)$  &   Volume of   sphere segment $|\{x=(x_1,x_2,x_3):x_1^2+x_2^2+x_3^2\leq r^2, x_1\leq t\}|=\frac{\pi}{3}(2r^3+3r^2t-t^3)$ \\\hline
			$\Pr\{\cdot\}$ & The probability of an event \\\hline
			{ $f(n)=O(g(n))$}   &  {  $f(n)\leq cg(n)$ for sufficiently large $n$ where $c>0$ }   \\\hline
			$f(n)=\Theta(g(n))$     &   $c_1g(n)\leq f(n)\leq c_2g(n)$ for sufficiently large $n$ where $0<c_1<c_2$     \\\hline
			$f(n)=o(g(n))$ &   $\lim\limits_{n\rightarrow\infty}\frac{f(n)}{g(n)}=0$ \\\hline
			$f(n)\sim g(n)$ &    $\lim\limits_{n\rightarrow\infty}\frac{f(n)}{g(n)}=1$   \\\hline
		\end{tabular}
	\end{center}
\end{table*}

\section{Preliminary results about regions  in $\mathbb{R}^3$}\label{sec:Preliminary-Results}

This section presents some fundamental characteristics of a region in $\mathbb{R}^3$.

\subsection{Spheres in  $\mathbb{R}^3$}

\begin{figure}[htbp]
	\centering
	\includegraphics[width=5.5cm]{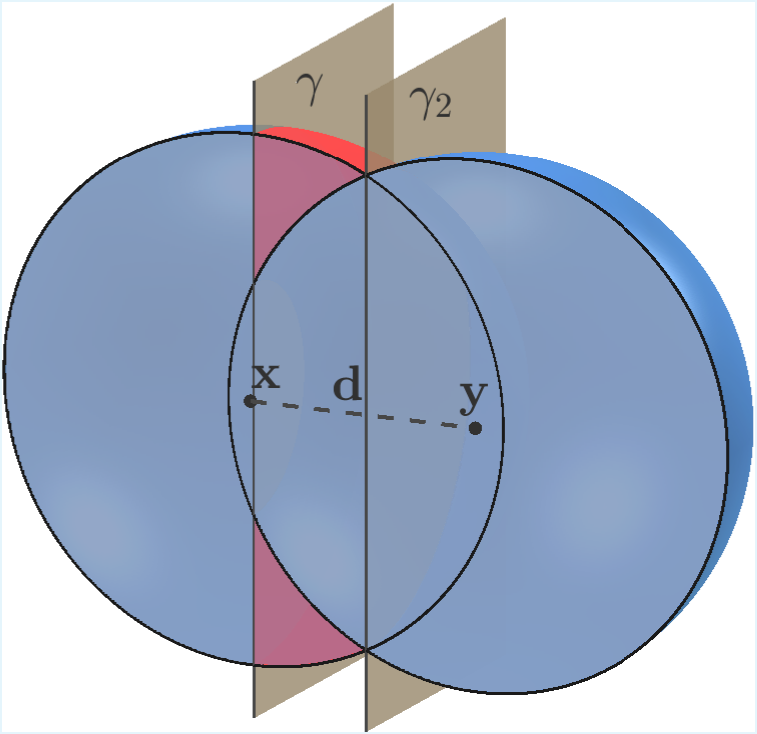}
	\caption{Illustration of Lemma~\ref{lem:shadow-low-bound} }
	\label{fig:Lemma3}
\end{figure}

\begin{lemma}\label{lem:shadow-low-bound}
Suppose  the centre  distance of two spheres $B(x,r)$ and $B(y,r)$  in $\mathbb{R}^3$
is  less than $r$, i.e.,
 $d=\mathrm{dist}(x,y)<r$.
	Plane $\gamma$ is perpendicular to line $xy$ and passes through point $x$. The hemisphere centred  at $x$,   cut by $\gamma$ and contains point $y$, is denoted by $B_{\mathrm{hemi}}(x,r)$, then the volume of $B_{\mathrm{hemi}}(x,r)\setminus B(y,r)$ (see red part in Figure~\ref{fig:Lemma3}) is
	$$V^{*}(d)=|B_{\mathrm{hemi}}(x,r)\setminus B(y,r)|=\frac{1}{4}\pi d^3.$$
\end{lemma}
\begin{proof} By simple calculation,
	$$
	V^{*}(d)=\int_0^{\frac{d}{2}}\pi\left[(r^2-t^2)-(r^2-(d-t)^2)\right]\mathrm{d}t=\frac{1}{4}\pi d^3.
	$$
\end{proof}

Throughout this article, for any $t\in[0,r]$, we define
\begin{equation}\label{eq:a(r,t)}
	a(r,t)=|\{x=(x_1,x_2,x_3):x_1^2+x_2^2+x_3^2\leq r^2, x_1\leq t\}|,
\end{equation}
the volume of the major  segment of a sphere cut by a plane.
Function $a(r,t)$  is frequently   used in this article, and    shortly denoted by $a(t)$  if there is no  confusion about $r$.
It is easy to see that
\begin{equation}\label{eq:a(r,t)-2}
	a(t)=\frac{\pi}{3}(r+t)^2(2r-t)=\frac{\pi}{3}(2r^3+3r^2t-t^3),
\end{equation}
\begin{equation}
a'(t)=\frac{\mathrm{d} a(r,t)}{\mathrm{d}t}=\pi(r^2-t^2), \quad  a''(t)=-2\pi t.
\end{equation}

The next lemma gives an estimate on   function   $a(t)$.
\begin{lemma}\label{lemma:a(t)}
	Let $r= \left( \frac{\log n+(\frac{3k}{2}-1)\log\log n+\xi}{ n\pi }\right)^{\frac13}$ for $k\geq 0$. Then
	$$
	n\int_0^{\frac r2}\frac{(na(t))^ke^{-na(t)}}{k!}\mathrm{d}t\sim
	\frac{e^{-\frac{2\xi}{3}}}{ \pi^{\frac{1}{3}}}\left(\frac{2}{3}\right)^k\frac{1}{k!}.
	$$
\end{lemma}

\begin{proof}
	Denote $C_k=\frac{3k-2}{2}$.
	Let $f(t)=na(t)$, then
	\begin{equation*}
		\begin{aligned}
			&n\int_{0}^{\frac r2} \frac{(f(t))^ke^{-f(t)}}{k!}\mathrm{d}t\\
			=&n\int_{0}^{\frac r2} \frac{(f(t))^ke^{-f(t)}}{k!} \frac{1}{na'(t)} \mathrm{d}f(t)\\
			=&   \int_{0}^{\frac r2} \frac{1}{a'(t)}\mathrm{d}\left(-e^{-f(t)} \sum_{i=0}^{k}\frac{(f(t))^i}{i!}  \right)\\
			=&-\frac{1}{a'(t)}e^{-f(t)}\sum_{i=0}^{k}\frac{(f(t))^i}{i!}|_{0}^{\frac r2}\\
			& -\int_{0}^{\frac r2} \frac{a''(t)}{(a'(t))^2}e^{-f(t)}\sum_{i=0}^{k}\frac{(f(t))^i}{i!}\mathrm{d}t.
		\end{aligned}
	\end{equation*}
	That is,
	\begin{equation*}
		\begin{aligned}
			& n\int_{0}^{\frac r2} \frac{(f(t))^ke^{-f(t)}}{k!}\mathrm{d}t\\
			= &  \frac{1}{a'(0)}e^{-f(0)}\sum_{i=0}^{k}\frac{(f(0))^i}{i!}
			-\frac{1}{a'(r)}e^{-f(r)}\sum_{i=0}^{k}\frac{(f(r))^i}{i!}   \\
			& -\int_{0}^{\frac r2} \frac{a''(t)}{(a'(t))^2}e^{-f(t)}\sum_{i=0}^{k}\frac{(f(t))^i}{i!}\mathrm{d}t. \\
		\end{aligned}
	\end{equation*}
	
	(i) First term.  Notice that
	$$
	a'(0)= \pi r^2= \pi^{1/3} \left(\frac{\log n +C_k\log\log n+\xi}{n}\right)^{\frac{2}{3}},
	$$
	$$
	f\left(0\right)=\frac{2\pi nr^3}{3}=\frac{2}{3}\left(\log n+C_k\log\log n+\xi\right),$$
	$$
	e^{-f(0)}=  \frac{1}{n^{\frac{2}{3}}}\frac{1}{(\log n)^{\frac{2C_k}{3}}} e^{-\frac{2\xi}{3}}.
	$$
	It is easy to see that
	\begin{equation*}
		\begin{aligned}
			&\frac{1}{a'\left(0\right)}e^{-f\left(0\right)}f\left(0\right)^k\\
			=& \frac{1}{\pi^{1/3} \left(\frac{\log n +C_k\log\log n+\xi}{n}\right)^{\frac{2}{3}}} \frac{1}{n^{\frac{2}{3}}}\frac{1}{(\log n)^{\frac{2C_k}{3}}} e^{-\frac{2\xi}{3}} \\
			&\cdot\left[ \frac{2}{3}\left(\log n+C_k\log\log n+\xi\right)\right]^k\\
			=& \frac{1}{\pi^{\frac{1}{3}}e^{\frac{2\xi}{3}}} \frac{\left[ \frac{2}{3}\left(\log n+C_k\log\log n+\xi\right)\right]^k}
			{\left( \log n +C_k\log\log n+\xi \right)^{\frac{2}{3}} (\log n)^{\frac{2C_k}{3}}} \\
			\sim &\frac{e^{-\frac{2\xi}{3}}}{ \pi^{\frac{1}{3}}}\left(\frac{2}{3}\right)^k,
		\end{aligned}
	\end{equation*}
	where $\frac{2C_k}{3}+\frac{2}{3}=k\geq 0$.  For $k>0$ and $0\leq i< k$,
	$$ -\frac{1}{a'\left(0\right)}e^{-f\left(0\right)}f\left(0\right)^i=o(1).$$
	Hence,
	\begin{equation}
		\begin{aligned}\label{eq:temp-sim-inLemma}
			\frac{1}{a'\left(0\right)}e^{-f\left(0\right)}\sum_{i=0}^k\frac{f\left(0\right)^i}{i!} \sim & \frac{1}{a'\left(0\right)}e^{-f\left(0\right)} \frac{f\left(0\right)^k}{k!} \\
			\sim &\frac{e^{-\frac{2\xi}{3}}}{ \pi^{\frac{1}{3}}}\left(\frac{2}{3}\right)^k\frac{1}{k!}.
		\end{aligned}
	\end{equation}

	(ii) Second term.  Noticing
	$$
	a'\left(\frac r2\right)=\frac{3\pi r^2}{4} = \Theta \left(\frac{\log n}{n}\right)^{\frac{2}{3}},
	$$
	$$
	f\left(\frac r2\right) =
	\frac{9}{8}\left(  \log n+ C_k \log\log n+\xi \right)\sim \Theta( \log n),
	$$
	$$
	e^{-f\left(\frac r2\right)}= \frac{1}{n^{\frac{9}{8}}}\frac{1}{(\log n)^{\frac{9C_k}{8}}} e^{-\frac{9\xi}{8}},
	$$
	so
	\begin{equation*}
		\begin{aligned}
			& \frac{1}{a'\left(\frac{r}{2}\right)}e^{-f\left(\frac r2\right)}\sum_{i=0}^k\frac{f\left(\frac r2\right)^i}{i!}\\
			= & \frac{ \Theta( \log n)^k}{ \Theta \left(\frac{\log n}{n}\right)^{\frac{2}{3}}}\frac{1}{n^{\frac{9}{8}}}\frac{1}{(\log n)^{\frac{9C_k}{8}}} e^{-\frac{9\xi}{8}} =o(1).
		\end{aligned}
	\end{equation*}

	(iii) Third term.
	When $t\leq \frac{r}{2}$,
	$$
	\left|\frac{a''(t)}{(a'(t))^3}\right|=\frac{2}{\pi^2}\frac{t}{(r^2-t^2)^3}\leq \frac{64}{27\pi^2}\frac{1}{r^5}=\frac{\Theta(1)}{r^5}.
	$$
	Then
	\begin{eqnarray*}
		&&\left|\int_{0}^{\frac r2} \frac{a''(t)}{(a'(t))^2}e^{-f(t)}\sum_{i=0}^{k}\frac{(f(t))^i}{i!}\mathrm{d}t\right|\\
		&=&\frac{1}{n}\left|\int_{0}^{\frac r2} \frac{a''(t)}{(a'(t))^3}e^{-f(t)}\sum_{i=0}^{k}\frac{(f(t))^i}{i!}\mathrm{d}f(t)\right|\\
		&\leq & \Theta(1)\frac{1}{nr^5} \left|\int_{0}^{\frac r2} e^{-f(t)}\sum_{i=0}^{k}\frac{(f(t))^i}{i!}\mathrm{d}f(t)\right|\\
		&\leq &\frac{ \Theta(1)}{nr^5}  \left|\int_{0}^{\frac r2} e^{-f(t)}\frac{(f(t))^k}{k!}\mathrm{d}f(t)\right|\\
		&= & \frac{\Theta(1)}{nr^5} \left|\int_{0}^{\frac r2} \mathrm{d}\left(-e^{-f(t)}\sum_{i=0}^{k}\frac{(f(t))^i}{i!}\right) \right|\\
		&\leq &\frac{ \Theta(1)}{nr^5}e^{-f(0)}\sum_{i=0}^{k}\frac{(f(0))^i}{i!}  \\
		&= & \frac{\Theta(1)}{nr^3}  =o(1),
	\end{eqnarray*}
	where $e^{-f(0)}\sum_{i=0}^{k}\frac{(f(0))^i}{i!} =\Theta(1) a'(0)= \Theta(1) r^2$  by (\ref{eq:temp-sim-inLemma}).

Therefore,
$$
n\int_{0}^{\frac r2} \frac{(f(t))^ke^{-f(t)}}{k!}dt
\sim \frac{e^{-\frac{2\xi}{3}}}{\pi^{\frac13}}\left(\frac{2}{3}\right)^{k}\frac{1}{k!}.
$$

\end{proof}

\subsection{Compact region $\Omega$ in $\mathbb{R}^3$}

\begin{figure}[htbp]
\begin{center}
\scalebox{0.45}[0.45]{
	\includegraphics{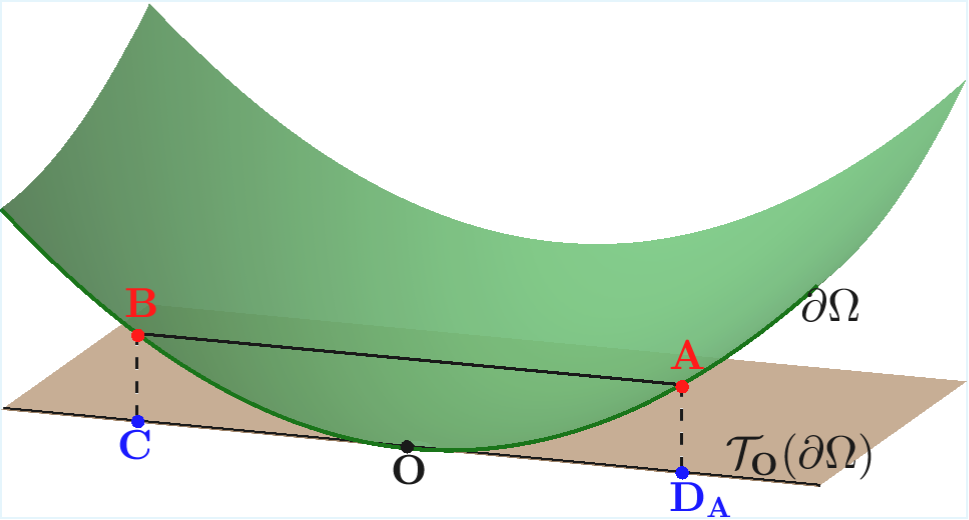}
}\caption{$\mathcal{T}_O(\partial\Omega)$ is tangent to $\partial \Omega$ at point $O$}\label{fig:Rectangle-ABCD}
\end{center}
\end{figure}

\begin{lemma}\label{lemma:G(Omega)}
 Suppose surface $\partial\Omega $ in $\mathbb{R}^3$ is twice continuously differentiable   in a neighborhood of  point $O\in\partial\Omega$, with tangent plane $\mathcal{T}_O(\partial\Omega)$ of   $\partial\Omega$ at point $O$.
   Let the projection of point $A\in\partial \Omega (A\neq O)$  onto the tangent plane     be $D_A\in\mathcal{T}_O(\partial\Omega)$, i.e.,
	line $AD_A \perp\mathcal{T}_O(\partial\Omega)$ at point $D_A$. Then
	\begin{equation}\label{eq:lem-G(Omega)-0}
		\|AD_A\|\leq   	G(\partial\Omega|_O)\|OD_A\|^2+o(\|OD_A\|^2),
	\end{equation}
	where
	$ 	G(\partial\Omega|_O)	$
	 is a constant only depending on point $O\in\partial\Omega$.
	This leads to
	that if $\|AD_A\|> (G(\partial\Omega|_O)+1)r^2$, then
	$\|OD_A\|> r, $
	as long as $r$ is sufficiently small.

	Furthermore, if  surface $\partial\Omega \subset \mathbb{R}^3$ is $C^2$ smooth, then there exists a universal   constant $G(\partial\Omega)$, which is  only determained by the curvature of $\partial\Omega$, such that  for any such point $O\in\partial\Omega$,
	\begin{equation}\label{eq:lem-G(Omega)-1}
	\|AD_A\|\leq   	G(\partial\Omega|)\|OD_A\|^2+o(\|OD_A\|^2),
\end{equation}	
and consequently, 	 if $\|AD_A\|> (G(\partial\Omega)+1)r^2$, then
$\|OD_A\|> r $ for sufficently small $r$.
\end{lemma}

\begin{proof}
	If $A\in\mathcal{T}_O(\partial\Omega)$, then $\|AD_A\|=0$. Otherwise, plane $\pi_{OAD_A}$ passing through three points $O,A,D_A$, is
	perpendicular to  tangent plane $\mathcal{T}_O(\partial\Omega)$, i.e.  $\pi_{OAD_A}\perp \mathcal{T}_O(\partial\Omega)$. Curve $\gamma =\pi_{OAD_A}\cap \partial\Omega$ passes points $O$ and $D_A$. According to the proof of Lemma~3 in~\cite{Ding2025-IEEETON},
	\begin{eqnarray}\label{eq:lemma3-Important}
		\|AD_A\|&=& \frac{1}{2} \mathrm{Curv}_O(\gamma)\cdot \|OD_A\|^2+o(\|OD_A\|^2),
	\end{eqnarray}
	where $\mathrm{Curv}_O(\gamma)$ is the curvature of curve $\gamma$ at point $O$.  The maximum and minimum principal curvatures of surface $\partial{\Omega}$ at point $O$ are denoted by $\mathrm{Curv_{max}}(\partial\Omega|_O)$
	and $\mathrm{Curv_{min}}(\partial\Omega|_O)$  respectively.  By the Euler's formula  and Cauchy inequality,
	\begin{equation*}
		\begin{aligned}
			&\mathrm{Curv}_O(\gamma)^2 \\
			=& \left[\mathrm{Curv_{max}}(\partial\Omega|_O)\cos^2\theta +\mathrm{Curv_{min}}(\partial\Omega|_O)\sin^2\theta \right]^2\\
			\leq & \left[\mathrm{Curv_{max}}(\partial\Omega|_O)^2 +\mathrm{Curv_{min}}(\partial\Omega|_O)^2\right]\left[\cos^4\theta+\sin^4\theta\right]\\
			\leq  &  \mathrm{Curv_{max}}(\partial\Omega|_O)^2 +\mathrm{Curv_{min}}(\partial\Omega|_O)^2.
		\end{aligned}
	\end{equation*}
	Let constant
		$$
	G(\partial\Omega|_O) =   \frac{  | \mathrm{Curv_{max}}(\partial\Omega|_O)|
		+ |\mathrm{Curv_{min}}(\partial\Omega|_O))|}{2}.
	$$
   Then $\frac{1}{2} \mathrm{Curv}_O(\gamma)\leq	G(\partial\Omega|_O) $, and inequality (\ref{eq:lem-G(Omega)-0}) holds.

\par Hereafter, for $C^2$ smooth surface $\partial\Omega \subset \mathbb{R}^3$, we  define
\begin{equation}\label{eq:G(Omega)}
	\begin{aligned}
		G(\partial\Omega) =& \sup_{O\in\partial\Omega}  	G(\partial\Omega|_O) \\
		=	&\frac{1}{2}\sup_{O\in\partial\Omega}  \left( | \mathrm{Curv_{max}}(\partial\Omega|_O)|
		+  |\mathrm{Curv_{min}}(\partial\Omega|_O))|\right).
	\end{aligned}
\end{equation}
Then  inequality (\ref{eq:lem-G(Omega)-1}) holds, and the proof is completed.
\end{proof}

This lemma states that if the distance from a point on  surface $\partial\Omega$ to the tangent plane at some surface point is less than $(G(\partial\Omega)+1)r^2$, then the distance from the projection of this  point onto the tangent plane to the tangent point is greater than $r$. We should point out that in this lemma,    $\partial\Omega$  can be a general surface, not necessarliy the boundary of a region.

\begin{figure}[htbp]
	\begin{center}
		\scalebox{0.5}{
			\includegraphics{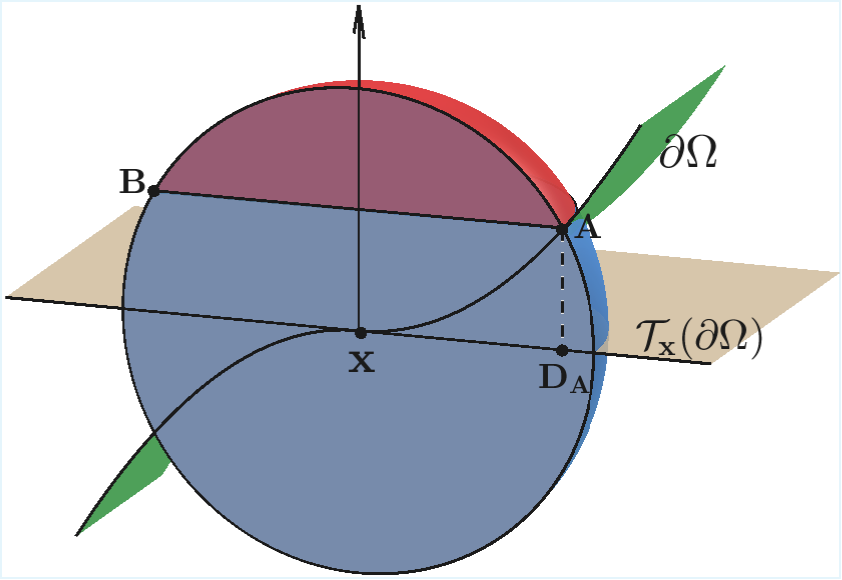}}
		\caption{Illustration of Lemma~\ref{lemma:C2-Boundary-volume}}\label{fig:lemma-C2-Boundary-volume}
	\end{center}
\end{figure}

\begin{lemma}\label{lemma:C2-Boundary-volume}
	(1). If  the boundary $\partial\Omega$ of  a compact region $\Omega\subset \mathbb{R}^3$ is
	 twice continuously differentiable   in a neighborhood of   point $x\in\partial\Omega$, then
	$$\liminf_{r\rightarrow 0}  \frac{\left|B(x,r)\cap \Omega\right|}{|B(x,r)|}\geq\frac{1}{2}+o(1), \quad r\rightarrow 0.$$
(2).	If the boundary $\partial\Omega$ of a compact region $\Omega\subset \mathbb{R}^3$ is $C^2$ smooth, then
	$$\liminf_{r\rightarrow 0} \inf_{x\in\Omega}\frac{\left|B(x,r)\cap \Omega\right|}{|B(x,r)|}\geq\frac{1}{2}+o(1), \quad r\rightarrow 0.$$

\end{lemma}
\begin{proof} We only prove (2). If $x$ is an interior point,  $B(x,r)\subset \Omega$  for sufficiently small $r$, then the conclusion holds. Let $x\in\partial\Omega$.  	Denote $\mathcal{T}_x(\partial\Omega)$ the tangent plane of $\partial\Omega$ at point $x$.
Suppose $A\in B(x,r)\cap\partial\Omega$, with its projection $D_A \in \mathcal{T}_x(\partial\Omega)$.  We may assume  $A$ is the farest point on  $B(x,r)\cap\partial\Omega$ from $\mathcal{T}_x(\partial\Omega)$.
 See Figure~\ref{fig:lemma-C2-Boundary-volume}. 	By Lemma~\ref{lemma:G(Omega)},
	\begin{equation*}
		\begin{aligned}
			\|AD_A\|\leq & G(\partial\Omega) \|xD_A\|^2+o(\|xD_A\|^2)\\
			<& (G(\partial\Omega)+1) \|xD_A\|^2\leq (G(\partial\Omega)+1) r^2.
		\end{aligned}
	\end{equation*}
Here constant $G(\partial\Omega)$ is given by (\ref{eq:G(Omega)}).
The minor spherical cap	of $B(x,r)$ cut by the plane that passes through point $A$ and is parallel to  $\mathcal{T}_x(\partial\Omega)$ has volume  $V=|B(x,r)|-a(\|AD_A\|)$, where function $a(\cdot)$ is defined by~(\ref{eq:a(r,t)}). See the red volume in Figure~\ref{fig:lemma-C2-Boundary-volume}.  The volume $V$ is bounded by
	\begin{equation*}
		\begin{aligned}
			V\geq &\frac{1}{2}|B(x,r)|-\pi r^2 \|AD_A\|  \\
			\geq & \frac{1}{2}|B(x,r)|-\pi(G(\partial\Omega)+1) r^4 =\left(\frac{1}{2}+o(1)\right) \frac{4}{3}\pi r^3,
		\end{aligned}
	\end{equation*}
leading to that
	$
	\frac{\left|B(x,r)\cap \Omega\right|}{|B(x,r)|}
	\geq \frac{V}{|B(x,r)|}
	\geq  \frac{1}{2}+o(1) ,    r\rightarrow 0.$
	
\end{proof}

\subsection{Bounds of $|B(x,r)\cap\Omega|$ under condition $C_{\mathrm{I}}$}\label{sec:B(x,r)-Omega-Condition-I}

In subsection, we will show that as long as a point $x\in\Omega$  is at an appropriate distance from the boundary $\partial\Omega$, then set $B(x,r)\cap\Omega$ contains the   major    segment of sphere $B(x,r)$   cut by a plane, and is contained in the union of this sphere segment and a circular cylinder of volume $\Theta(r^4)$. As a consequent,
the lower and upper bounds of $|B(x,r)\cap\Omega|$ will be  derived straightforwardly.

Suppose compact region $\Omega \subset \mathbb{R}^3$ has $C^2$ smooth boundary $\partial \Omega$ in this subsetion.
For any point $x\in \Omega$ satisfying
$$
(G(\partial\Omega)+1)r^2< \mathrm{dist}(x,\partial \Omega)<r,
$$
let $W_x\in\partial \Omega$ be the nearest point to the given $x\in\Omega$, so
$
\|xW_x\|=\mathrm{dist}(x,\partial \Omega)<r.
$
The tangent plane  $ \mathcal{T}_{W_x}(\partial\Omega)$ of $\partial\Omega$ at point $W_x$ is  shortly denoted by   $ \mathcal{T}_{W_x}$.

Let $\gamma_P(\mathcal{T}_{W_x} )$ denote the plane which passes through point $P$ and is parallel to  $\mathcal{T}_{W_x}$.
In particular,   plane $\gamma_x(\mathcal{T}_{W_x})$ contains point $x$  and
 $\gamma_x(\mathcal{T}_{W_x})\pll \mathcal{T}_{W_x}$.

Noticing $\|xW_x\|=\mathrm{dist}(x,\partial \Omega)>(G(\partial\Omega)+1)r^2$,  so by  Lemma~\ref{lemma:G(Omega)},
$$
\mathrm{dist}(x, \partial \Omega \cap \gamma_x(\mathcal{T}_{W_x}))>r.
$$
This implies that $\forall A\in\partial\Omega \cap \gamma_x(\mathcal{T}_{W_x})$,  $\|xA\|>r.$

\begin{figure}[htbp]
	\begin{center}
		\scalebox{0.43}{
			   \includegraphics{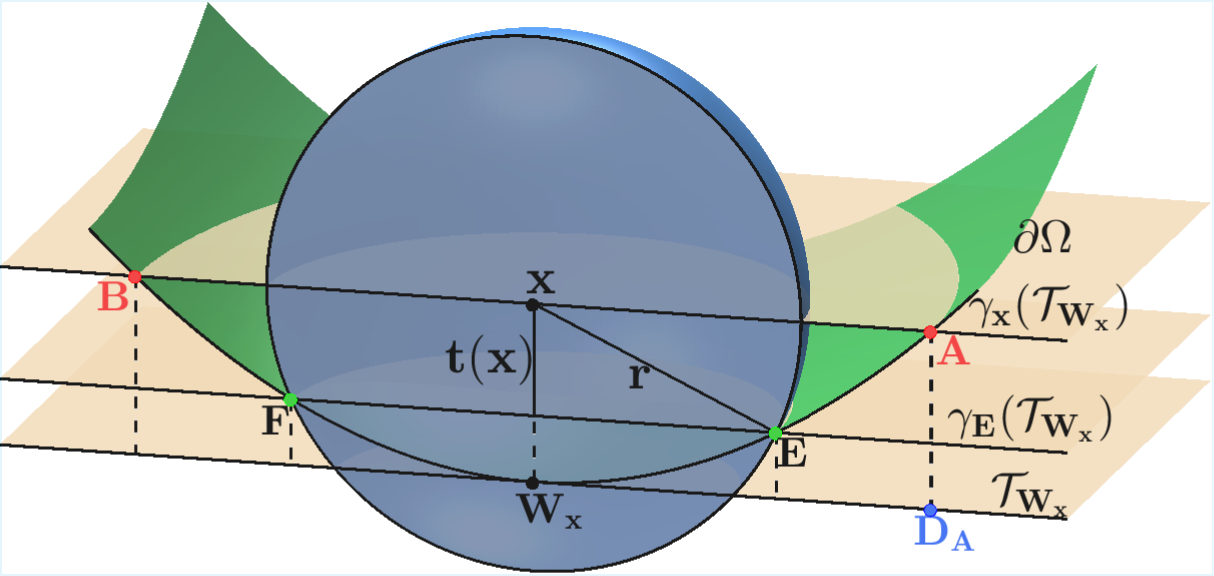}}
				\caption{$t(x)$ is nonnegative and major spherical segment contained in $\Omega$}\label{positive-at)}
	\end{center}
\end{figure}

The distance
between point $x$ and   point $y\in\partial \Omega$ is a continuous
function of $y$, which is due to the continuous
boundary $\partial \Omega$. By $\|xW_x\|<r$ and
$\|xA\|>r$, we know that there exits a point $E'\in \partial \Omega$ such that $\|xE'\|=r$, i.e., $E'\in\partial\Omega\cap\partial B(x,r)$.
Let
\begin{equation*}
	\begin{aligned}
		&t(x,r)
		= \inf\{\mathrm{dist}(x,\gamma_{E'}(\mathcal{T}_{W_x}))\mid E'\in\partial\Omega\cap\partial B(x,r)\},
	\end{aligned}
\end{equation*}
where is  shortly denoted by $t(x)$. It is easy to see  $t(x)>0$.
 See Figure~\ref{positive-at)} and~\ref{fig:Bounds-B(x,r)Omega}.

Now we give a lower bound for $|B(x,r)\cap \Omega|$. Because $\partial\Omega\cap\partial B(x,r)$ is compact, we assume $E \in\partial\Omega\cap\partial B(x,r)$ such that
\begin{equation*}
	\begin{aligned}
		  t(x)=&\mathrm{dist}(x,\gamma_{E}(\mathcal{T}_{W_x}))\\
         =&\inf\{\mathrm{dist}(x,\gamma_{E'}(\mathcal{T}_{W_x}))\mid E'\in\partial\Omega\cap\partial B(x,r) \}.
	\end{aligned}
\end{equation*}
That is,  plane  $\gamma_{E}(\mathcal{T}_{W_x})\pll \mathcal{T}_{W_x}$,   passes through point $E$ with $\mathrm{dist}(x, \gamma_{E}(\mathcal{T}_{W_x}))=t(x)$.
The   major    segment of sphere $B(x,r)$ cut by $\gamma_{E}(\mathcal{T}_{W_x})$ is denoted by $B(x,r)|_{\gamma_{E}(\mathcal{T}_{W_x})}$, then
\begin{equation}\label{eq:B(x,r)-Omega-subset}
B(x,r)|_{\gamma_{E}(\mathcal{T}_{W_x})}\subset B(x,r)\cap\Omega,
\end{equation}
and thus
\begin{equation}\label{eq:volume-major-sphere segment}
	\begin{aligned}
		|B(x,r)\cap\Omega|\geq& | B(x,r)|_{\gamma_{E}(\mathcal{T}_{W_x})}|\\
		=&a(t(x))=  \frac{\pi}{3}(2r^3+3r^2t(x)-t(x)^3),
	\end{aligned}
\end{equation}
where $a(t(x))$ is specified by~(\ref{eq:a(r,t)}).

\par
\begin{figure}[htbp]
	\begin{center}
		\scalebox{0.42}{
			   \includegraphics{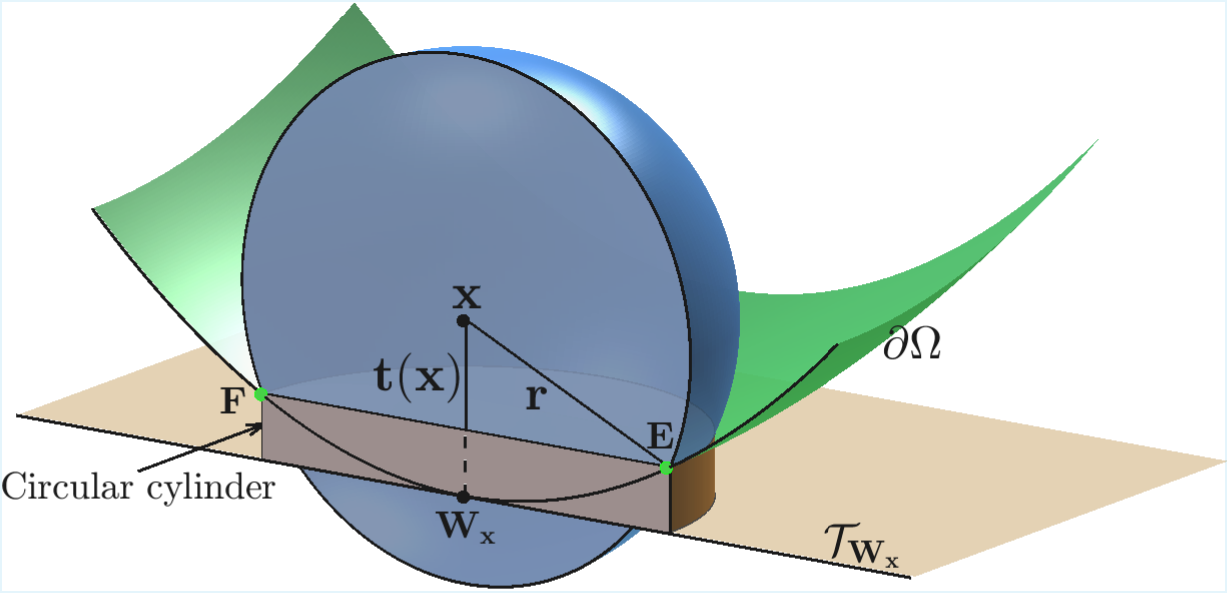}}
				\caption{Upper bound for $|B(x,r)\cap\Omega|$}\label{fig:UpperBound-Ball-Omega}
	\end{center}
\end{figure}

 Now we consider   upper bounds of  $|B(x,r)\cap\Omega|$.
See Figure~\ref{fig:UpperBound-Ball-Omega} and~\ref{fig:Bounds-B(x,r)Omega-2}.
Notice that curve $\partial B(x,r)\cap \gamma_{E}(\mathcal{T}_{W_x})$ is a circle on plane
$\gamma_{E}(\mathcal{T}_{W_x})$. The corresponding disk enclosed by this circle is denoted by  $\mathrm{disk}(\partial B(x,r)\cap \gamma_{E}(\mathcal{T}_{W_x}))$.
There are two cases to discuss.

\par (1).  Assume all points in $\partial\Omega\cap \partial B(x,r)$  lie on the same side of tangent plane $\mathcal{T}_{W_x}$  as point $x$.  Consider the projection of   $\mathrm{disk}(\partial B(x,r)\cap \gamma_{E}(\mathcal{T}_{W_x}))$   onto $\mathcal{T}_{W_x}$, which is denoted by    $\mathrm{Proj}\left( \mathrm{disk}(\partial B(x,r)\cap \gamma_{E}(\mathcal{T}_{W_x}))\right)|_{\mathcal{T}_{W_x}}$.
By $ \mathrm{Cyc}\left(\mathrm{disk}(\partial B(x,r)\cap \gamma_{E}(\mathcal{T}_{W_x})|_{\gamma_{E}(\mathcal{T}_{W_x}), \mathcal{T}_{W_x} }\right)$, we denote the circular cylinder formed by
$$
\mathrm{disk}(\partial B(x,r)\cap \gamma_{E}(\mathcal{T}_{W_x})) \subset\gamma_{E}(\mathcal{T}_{W_x})
$$
and its projection
$$
\mathrm{Proj}\left( \mathrm{disk}(\partial B(x,r)\cap \gamma_{E}(\mathcal{T}_{W_x}))\right)|_{\mathcal{T}_{W_x}}\subset \mathcal{T}_{W_x}.
$$
The volume of this circular cylinder is
\begin{equation*}
	\begin{aligned}
         & |\mathrm{disk}(\partial B(x,r)\cap \gamma_{E}(\mathcal{T}_{W_x}))|   \mathrm{dist}(\gamma_{E}(\mathcal{T}_{W_x}),\mathcal{T}_{W_x})\\
		  =&\pi (r^2-t(x)^2)\mathrm{dist}(\gamma_{E}(\mathcal{T}_{W_x}),\mathcal{T}_{W_x})\\
		\leq &\pi (r^2-t(x)^2) (G(\partial\Omega)+o(1))r^2\\
		\leq & \pi (G(\partial\Omega)+o(1))r^4.
	\end{aligned}
\end{equation*}
Here  by Lemma~\ref{lemma:G(Omega)},
\begin{equation}\label{eq:dist-planes}
	\begin{aligned}
		&\mathrm{dist}(\gamma_{E}(\mathcal{T}_{W_x}),\mathcal{T}_{W_x}(\partial\Omega))\\
		\leq & G(\partial\Omega) (\|xE\|^2) +o(\|xE\|^2)) =(G(\partial\Omega)+o(1))r^2.
	\end{aligned}
\end{equation}
Notice that
\begin{equation*}
	\begin{aligned}
	 B(x,r)\cap\Omega
		\subset &  B(x,r)|_{\gamma_{E}(\mathcal{T}_{W_x})} \bigcup \\		
		& \mathrm{Cyc}\left(\mathrm{disk}(\partial B(x,r)\cap \gamma_{E}(\mathcal{T}_{W_x})|_{\gamma_{E}(\mathcal{T}_{W_x}), \mathcal{T}_{W_x} }\right).
	\end{aligned}
\end{equation*}
As a result,
\begin{equation}\label{eq:temp-upper-bound-1}
	\begin{aligned}
		&|B(x,r)\cap \Omega|\\
		\leq& | B(x,r)|_{\gamma_{E}(\mathcal{T}_{W_x})} |\\
		+ &|\mathrm{Cyc}\left(\mathrm{disk}(\partial B(x,r)\cap \gamma_{E}(\mathcal{T}_{W_x})|_{\gamma_{E}(\mathcal{T}_{W_x}), \mathcal{T}_{W_x} }\right)|\\
		\leq & a(t(x))+\pi (G(\partial\Omega)+o(1))r^4. \\
	\end{aligned}
\end{equation}

\begin{figure}[htbp]
	\begin{center}
		\scalebox{0.32}[0.32]{
			\includegraphics{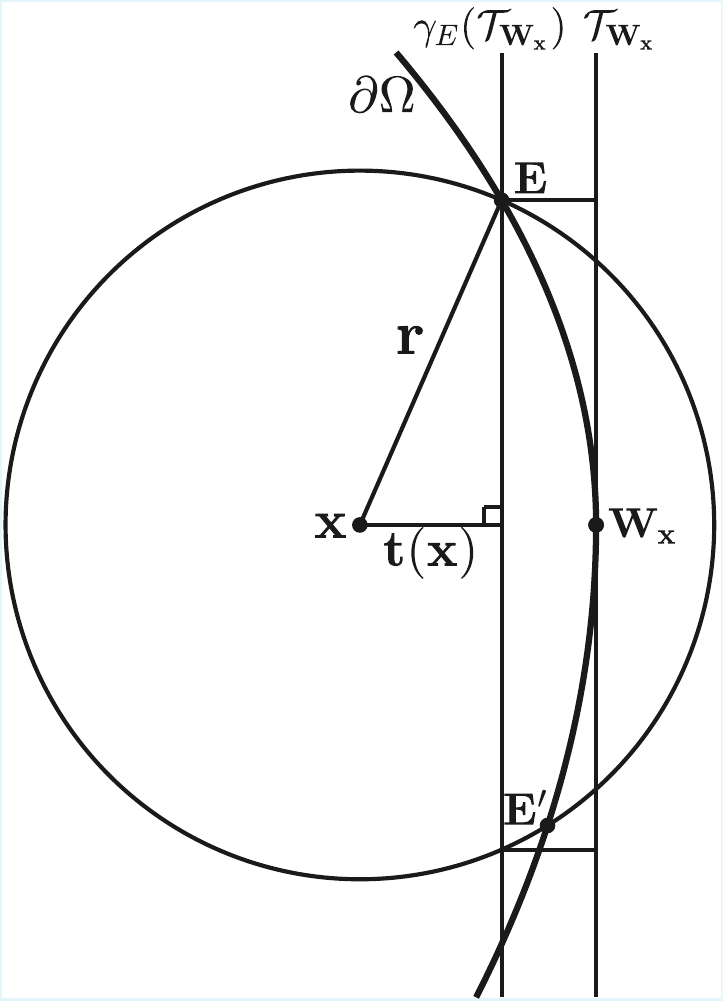}       }
		\caption{Lower and Upper bounds of $|B(x,r)\cap
			\Omega|$: a 2D illustration (1)}\label{fig:Bounds-B(x,r)Omega}
	\end{center}
\end{figure}

\begin{figure}[htbp]
	\begin{center}
		\scalebox{0.32}[0.32]{
			\includegraphics{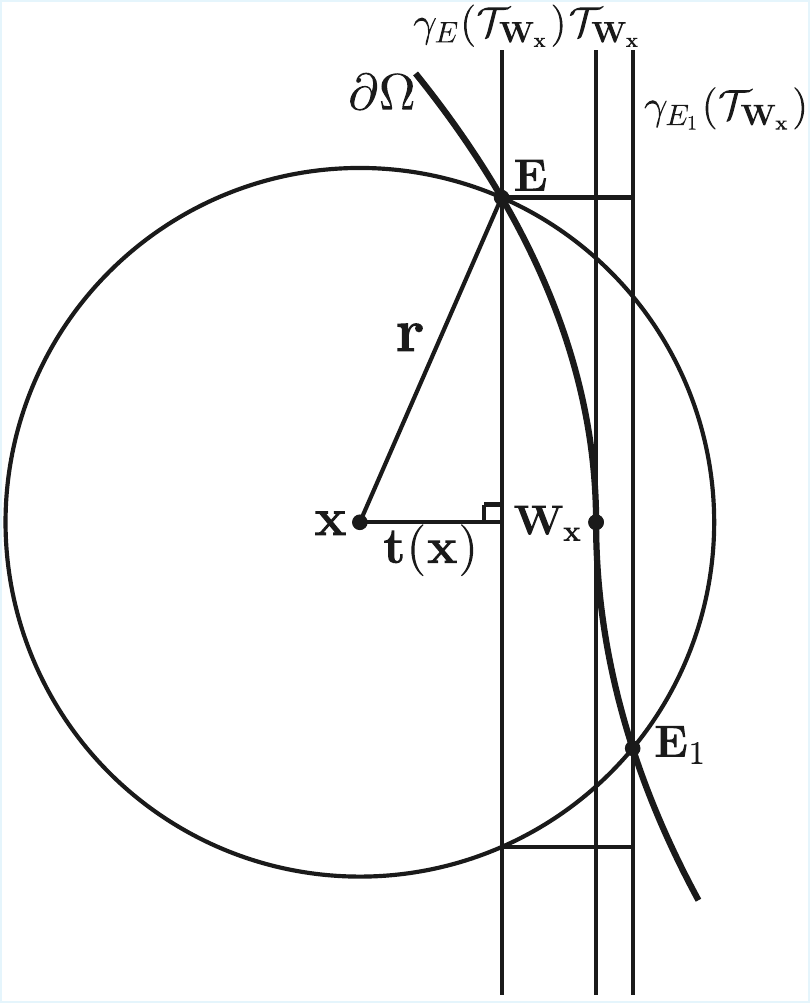}}
		\caption{Lower and Upper bounds for $|B(x,r)\cap
			\Omega|$: a 2D illustration (2)}\label{fig:Bounds-B(x,r)Omega-2}
	\end{center}
\end{figure}

\par (2). Suppose that there exit  some point   $ E''\in\partial\Omega\cap\Omega$ that lie on the opposite side of tangent plane $\mathcal{T}_{W_x}$ from $x$. We turn to consider the planes $\gamma_{E''}(\mathcal{T}_{W_x})\pll \mathcal{T}_{W_x}$.
 By similar discussions, we may assume   point $E_1\in\partial\Omega\cap\partial B(x,r)$, such that
plane $\gamma_{E_1}(\mathcal{T}_{W_x})$ to achieve the following maximum distance from point $x$:
\begin{equation*}
	\begin{aligned}
		&\mathrm{dist}(x,\gamma_{E_1}(\mathcal{T}_{W_x}))\\
= &\sup\{\mathrm{dist}(x,\gamma_{E''}(\mathcal{T}_{W_x})) \mid E''\in\partial\Omega\cap\partial B(x,r) \}.
	\end{aligned}
\end{equation*}
Similarly to (\ref{eq:dist-planes}), we have
$$
\mathrm{dist}(\gamma_{E_1}(\mathcal{T}_{W_x}),\mathcal{T}_{W_x})
 \leq  \pi  (G(\partial\Omega)+o(1))r^2.
$$
Hence,
\begin{equation*}
	\begin{aligned}
		&\mathrm{dist}(\gamma_{E}(\mathcal{T}_{W_x},\gamma_{E_1}(\mathcal{T}_{W_x})\\
		=&\mathrm{dist}(\gamma_{E}(\mathcal{T}_{W_x}),\mathcal{T}_{W_x})
		 +\mathrm{dist}(\gamma_{E_1}(\mathcal{T}_{W_x}),\mathcal{T}_{W_x})\\
		\leq &2\pi(G(\partial\Omega)+o(1))r^2.
	\end{aligned}
\end{equation*}
Similarly, we can obtain a circular cylinder $\mathrm{Cyc}\left(\mathrm{disk}(\partial B(x,r)\cap \gamma_{E}(\mathcal{T}_{W_x})|_{\gamma_{E}(\mathcal{T}_{W_x}); \gamma_{E_1}(\mathcal{T}_{W_x})}  \right)$, which is formed by  $\mathrm{disk}(\partial B(x,r)\cap \gamma_{E,W}) \subset\gamma_{E}(\mathcal{T}_{W_x})$ and its projection
on  $\gamma_{E_1}(\mathcal{T}_{W_x})$. Its  volume  is
\begin{equation*}
	\begin{aligned}
		&\left|\mathrm{Cyc}\left(\mathrm{disk}(\partial B(x,r)\cap \gamma_{E}(\mathcal{T}_{W_x})|_{\gamma_{E}(\mathcal{T}_{W_x}); \gamma_{E_1}(\mathcal{T}_{W_x})}  \right)\right|\\
		\leq & 2\pi(G(\partial\Omega)+o(1))r^2 \pi (r^2-t(x)^2).
	\end{aligned}
\end{equation*}
Because
\begin{equation*}
	\begin{aligned}
		&B(x,r)\cap\Omega \\
\subset&  B(x,r)|_{\gamma_{E}(\mathcal{T}_{W_x})} \bigcup  \\
 &\mathrm{Cyc}\left(\mathrm{disk}(\partial B(x,r)\cap \gamma_{E}(\mathcal{T}_{W_x})|_{\gamma_{E}(\mathcal{T}_{W_x}); \gamma_{E_1}(\mathcal{T}_{W_x})}  \right),
	\end{aligned}
\end{equation*}
we similarly obtain
\begin{equation}\label{eq:temp-upper-bound-2}
	\begin{aligned}
		&|B(x,r)\cap \Omega| \\
		\leq& | B(x,r)|_{\gamma_{E}(\mathcal{T}_{W_x})} |\\
		+&\left|\mathrm{Cyc}\left(\mathrm{disk}(\partial B(x,r)\cap \gamma_{E}(\mathcal{T}_{W_x}))|_{\gamma_{E}(\mathcal{T}_{W_x}); \gamma_{E_1}(\mathcal{T}_{W_x})}  \right)\right|\\
		\leq & a(t(x))+ 2\pi(G(\partial\Omega)+o(1))r^4. \\
	\end{aligned}
\end{equation}
Combine inequalities  (\ref{eq:temp-upper-bound-1}) and (\ref{eq:temp-upper-bound-2}),
we have obtained the following upper bound holding for all these two cases,
\begin{equation}\label{eq:temp-upper-bound}
	\begin{aligned}
		|B(x,r)\cap \Omega| \leq& a(t(x,r))+2\pi(G(\partial\Omega)+o(1))r^4\\
		=&  a(t(x,r))+\Theta(r^4)\\
		=&  \frac{\pi}{3}(2r^3+3r^2t(x)-t(x)^3)+\Theta(r^4).
	\end{aligned}
\end{equation}
Due to (\ref{eq:volume-major-sphere segment}) and (\ref{eq:temp-upper-bound}),
we have proved the following lemma.
\begin{lemma}\label{lemma:Bounds-B(x,r)-Omega}
	Suppose compact region $\Omega \subset \mathbb{R}^3$ has $C^2$ smooth boundary $\partial \Omega$.
	For any point $x\in \Omega$ satisfying $ (G(\partial\Omega)+1)r^2< \mathrm{dist}(x,\partial \Omega)<r$, where constant $G(\partial\Omega)$ is given by Lemma~\ref{lemma:G(Omega)}, let $W_x\in\partial \Omega$ be its
	 nearest point on boundary $\partial\Omega$ to $x$,
	  with $ \|xW_x\|=\mathrm{dist}(x,\partial \Omega)$. The tangent plane
	of $\partial\Omega$ at $W_x$  is   $\mathcal{T}_{W_x}(\partial\Omega)$  or shortly $\mathcal{T}_{W_x}$.
	Denote  by $\gamma_{E'}(\mathcal{T}_{W_x})$   the plane  that passes through point $E'$ and is parallel to   $\mathcal{T}_{W_x} $.
	Define
	\begin{equation}\label{eq:t(x,r)}
		 t(x,r)=\inf\{\mathrm{dist}(x,\gamma_{E'}(\mathcal{T}_{W_x}))\mid E'\in\partial\Omega\cap\partial B(x,r)   \}.
	\end{equation}
Then
	$|B(x,r)\cap\Omega|$ has the following bounds:
	$$
	    a(t(x,r))\leq  	|B(x,r)\cap\Omega|\leq   a(t(x,r))+\Theta(r^4),
	$$
where $	
  a(t(x,r))= \frac{\pi}{3}(2r^3+3r^2t(x,r)-t(x,r)^3),
	$
specified by  (\ref{eq:a(r,t)}).
\end{lemma}

\subsection{Bounds of $|B(x,r)\cap\Omega|$ under condition $C_{\mathrm{II}}$}\label{sec:B(x,r)-Omega-bound-sigma}

Lemma~\ref{lemma:Bounds-B(x,r)-Omega} essentially states that, provided
$x\in\Omega$ is propoerly  close to the boundary $\partial\Omega$
and at the nearest point on
 $\partial\Omega$ to $x$ the function of $\partial\Omega$ is twice continuously differentiable in a neighborhood, the set
$B(x,r)\cap \Omega$ can be bounded  using a spherical cap and a cylinder.
Indeed,  Lemma~\ref{lemma:Bounds-B(x,r)-Omega} is a local conclusion.
In the proof of Lemma~\ref{lemma:Bounds-B(x,r)-Omega}, the  sphere $B(x,r)$ under consideration is local, and the smoothness condition is applied locally.

 \par This subsection aims to establish   similar results under some
non-global smooth conditions like  condition $C_{\mathrm{II}}$.
The boundary $\partial\Omega$ of region $\Omega\subset \mathbb{R}^3$ is  $C^2$ smooth everywhere, except for a curve $\sigma\subset\partial\Omega$, where  $\sigma$ has zero area but finite length, i.e.,
  $\mathrm{Len}(\sigma)<\infty$ and    $\mathrm{Are}(\sigma)=0$.
We   use notation $G(\partial\Omega\backslash\sigma)$ to represent the   curvature constant of surface $\partial\Omega\backslash\sigma$, which can be similarly derived for surface $\partial\Omega\backslash\sigma$, according to the proof of  Lemma~\ref{lemma:G(Omega)}.   Due to  the assumption  of uniformly bounded maxisum and minimum principle curvatures  in condition $C_{\mathrm{II}}$,   $G(\partial\Omega\backslash\sigma)<\infty$.

Further, we define some tubes containing curve $\sigma$:
\begin{equation}\label{eq:tubes-sigma}
	\begin{aligned}
		\mathrm{Tube}_{2r}(\sigma)=&\{x\in\Omega: \mathrm{dist}(x,\sigma)\leq 2r \},\\
		\mathrm{Tube}_{3r}(\sigma)=&\{x\in\Omega: \mathrm{dist}(x,\sigma)\leq3r\}. 	
	\end{aligned}
\end{equation}
Consider any point  $x\in\Omega\backslash \mathrm{Tube}_{3r}(\sigma)$, satisfying
$$
(G(\partial\Omega\backslash\sigma)+1)r^2
< \mathrm{dist}(x,\partial \Omega)<r.
$$
We have
$$
(G(\partial\Omega\backslash\sigma)+1)r^2
< \mathrm{dist}(x,\partial \Omega)<r<3r<\mathrm{dist}(x,\sigma).
$$
Let $W_x\in\partial\Omega$ be the nearest point on $\partial \Omega$  to   $x$, with $\|xW_x\|=\mathrm{dist}(x,\partial\Omega)$.
Notice that
$$
|B(x,r)\cap\Omega|\subset \Omega\backslash \mathrm{Tube}_{2r}(\sigma),
$$
and $ W_x\in \partial\Omega\backslash  \mathrm{Tube}_{2r}(\sigma)$, which is due to
$$
\mathrm{dist}(W_x,\sigma)>\mathrm{dist}(x,\sigma)-\mathrm{dist}(x,W_x)>2r.
$$
That is, both $x$ and $W_x$ are sufficiently far from set $\sigma$, and hence  far from the singularity part of the boundary. Then, repeat the proof procedure of  Lemma~\ref{lemma:Bounds-B(x,r)-Omega}, we can similarly obtain the following conclusion:

\begin{corollary}\label{corollary:Bounds-B(x,r)-Omega-sigma}
	Suppose the boundary $\partial \Omega$ of a compact region $\Omega \subset \mathbb{R}^3$  is  $C^2$ smooth everywhere, except for a small part $\sigma\subset\partial\Omega$, where  $\mathrm{Len}(\sigma)<\infty$ and    $\mathrm{Are}(\sigma)=0$. The  curvature constant  (given according to  Lemma~\ref{lemma:G(Omega)}) of    surface $\partial\Omega\backslash\sigma$ is finite, i.e.	 $G(\partial\Omega\backslash\sigma)<\infty$.
	For any point $x\in  \Omega\backslash \mathrm{Tube}_{3r}(\sigma)$ satisfying $ (G(\partial\Omega\backslash\sigma)+1)r^2< \mathrm{dist}(x,\partial \Omega)<r$,
	$|B(x,r)\cap\Omega|$ is bounded by
	$$
	a(t(x,r))\leq  	|B(x,r)\cap\Omega|\leq   a(t(x,r))+\Theta(r^4),
	$$
	where $t(x,r)$ and $
	a(t(x,r))= \frac{\pi}{3}(2r^3+3r^2t(x,r)-t(x,r)^3),
	$ are defined similarly to   Lemma~\ref{lemma:Bounds-B(x,r)-Omega}.
\end{corollary}

\section{Proof of Proposition~\ref{pro:Explicit-Form} under condition $C_{\mathrm{I}}$}\label{sec:Pro-1}

This section is to give a proof of  Proposition~\ref{pro:Explicit-Form} under condition $C_{\mathrm{I}}$,  i.e., $ n\int_{\Omega} \psi^k_{n,r}(x)\mathrm{d}x\sim e^{-c}$,  where
the boundary $\partial \Omega$   is $C^2$ smooth.

Define
$$
\Omega(0)=\left\{x\in\Omega:\mathrm{dist}(x,\partial\Omega)\geq r\right\}, $$
$$
\Omega(2)=\left\{x\in\Omega:\mathrm{dist}(x,\partial\Omega)\leq
(G(\partial\Omega)+1)r^2\right\},
$$
where $r$ is considered sufficiently small  so that
$(G(\partial\Omega)+1)r^2<r$, resulting in $\Omega(0)\cap \Omega(2)=\emptyset$. Thus, let
\begin{equation*}
	\begin{split}
		\Omega(1)=&\Omega\backslash\left(\Omega(0)\cup\Omega(2)\right)\\
		=&\left\{x\in\Omega:
		(G(\partial\Omega)+1)r^2 < \mathrm{dist}(x,\partial\Omega)< r^2
		\right\}.
	\end{split}
\end{equation*}
Here constant $G(\partial\Omega)$ is determined in  Lemma~\ref{lemma:G(Omega)}.
Therefore,
\begin{equation*}
	\begin{split}
		&n\int_{\Omega}\psi^k_{n,r}(x) \mathrm{d}x= \left(\int_{\Omega(0)}+\int_{\Omega(2)}+\int_{\Omega(1)}\right)n\psi^k_{n,r}(x) \mathrm{d}x.
	\end{split}
\end{equation*}
The integrations   will be estimated
in the following Claim~\ref{pro:Omega(0)}-\ref{pro:Omega(1,1)}.

\begin{claim}\label{pro:Omega(0)}
	$n\int_{\Omega(0)} \psi^k_{n,r}(x)\mathrm{d}x= o\left(1\right). $
\end{claim}
\begin{proof}
	$\forall x\in \Omega(0), |B(x,r)\cap\Omega|=\frac{4}{3}\pi r^3$.  Notice that $|\Omega(0)|<1$ and by Remark~\ref{rem:remark3},
	\begin{eqnarray*}
		n\int_{\Omega(0)} \psi^k_{n,r}(x)\mathrm{d}x &=&\frac n{k!}\left(\frac{4n\pi r^3}{3}\right)^k e^{-\frac{4n\pi r^3}{3}}|\Omega(0)|\\
		& \leq & \frac n{k!}\left(\frac{4n\pi r^3}{3}\right)^k e^{-\frac{4n\pi r^3}{3}} = o\left(1\right).
	\end{eqnarray*}
\end{proof}

\begin{claim}\label{pro:Omega(2)}
	$
	n\int_{\Omega(2)} \psi^k_{n,r}(x)\mathrm{d}x=o\left(1\right).
	$
\end{claim}
\begin{proof} Notice that
	$$|\Omega(2)|\leq \mathrm{Area}(\partial\Omega)\times (G(\partial\Omega)+1)r^2=\Theta(1)r^2
	=\Theta \left(\frac{\log n}{n}\right)^\frac{2}{3}.$$
%
Here
	$n\pi r^3=\Theta\left( \log n\right)$, see Remark~\ref{rem:remark3}.
\begin{equation*}
\begin{aligned}
		n\int_{\Omega(2)} \psi^k_{n,r}(x)\mathrm{d}x  \leq &\frac n{k!}(n\pi r^3)^ke^{-\frac{2}{3}n\pi r^3}  |\Omega(2)|\\
		  =&\Theta(1)\frac{\left( \log n\right)^{k+\frac{2}{3}}}{n^{\frac{1}{3}}}=o(1).
\end{aligned}
\end{equation*}
\end{proof}

For any point $x\in \Omega(1)$,
$
(G(\partial\Omega)+1)r^2< \mathrm{dist}(x,\partial \Omega)<r.
$
According to the  lower and upper bounds of $|B(x,r)\cap \Omega|$ given in Lemma~\ref{lemma:Bounds-B(x,r)-Omega}:
$$
a(t(x,r))\leq  	|B(x,r)\cap\Omega|\leq   a(t(x,r))+\Theta(r^4),
$$
we have the
following estimates:
\begin{equation}\label{eq:local-1}
	\begin{split}
		&(1+o(1))(na(t(x)))^ke^{-na(t(x))}\\	
		= & e^{-n \Theta(r^4)}(na(t(x)))^ke^{-na(t(x))} \\
		\leq & (na(t(x)))^ke^{-n(a(t(x))+\Theta(r^4))}\\		
		\leq &(n|B(x,r)\cap \Omega|)^ke^{-n|B(x,r)\cap \Omega|}\\
		\leq& (1+o(1))^k(na(t(x)))^ke^{-na(t(x))},
	\end{split}
\end{equation}
or equivalently,
\begin{equation}\label{eq:estimation-psi(n)}
	(n|B(x,r)\cap \Omega|)^ke^{-n|B(x,r)\cap \Omega|} \sim   (na(t(x)))^ke^{-na(t(x))}.
\end{equation}

Set
$$
\Omega(1,1)=\left\{x\in \Omega(1):t(x)\leq \frac r2\right\},
$$
$$
\Omega(1,2)=\Omega(1)\setminus\Omega(1,1),
$$
where $t(x)=t(x,r)$ is defined in Lemma~\ref{lemma:Bounds-B(x,r)-Omega}. Then
\begin{equation*}
	n\int_{\Omega(1)}\psi^k_{n,r}(x)\mathrm{d}x=  \int_{\Omega(1,1)}n\psi^k_{n,r}(x)\mathrm{d}x+ \int_{\Omega(1,2)}n\psi^k_{n,r}(x)\mathrm{d}x.
\end{equation*}
In following, we will specify $n\int_{\Omega(1,2)}\psi^k_{n,r}(x)\mathrm{d}x$ in
Claim~\ref{pro:Omega(1,2)}, and then  determine $n\int_{\Omega(1,1)}\psi^k_{n,r}(x)\mathrm{d}x$ in
Claim~\ref{pro:Omega(1,1)}.

\begin{claim}\label{pro:Omega(1,2)}
	$
	n\int_{\Omega(1,2)} \psi^k_{n,r}(x)\mathrm{d}x=o(1).
	$
\end{claim}

\begin{proof}
	Notice that volume
	$|\Omega(1,2)|\leq \mathrm{Area}(\partial\Omega)r\sim \mathrm{Area}(\partial\Omega) \left(\frac{\log n}{n}\right)^{\frac{1}{3}}$.
	Because $t(x)\geq \frac{r}{2}$ for any $x\in \Omega(1,2)$,  so $
	a(t(x,r))\geq \frac{2}{3}\pi r^3+\frac{1}{2}\frac{2}{3}\pi r^3=\pi r^3
	$. Then by formula~(\ref{eq:local-1}),
	\begin{eqnarray*}
		&&n\int_{\Omega(1,2)} \psi^k_{n,r}(x)\mathrm{d}x\\
		&=&\frac n{k!}\int_{\Omega(1,2)}\left(n|B(x,r_n)\cap\Omega|\right)^k e^{-n|B(x,r_n)\cap\Omega|}\mathrm{d}x\\
		&\leq&\frac n{k!}\left[(1+o(1))^k(na(t(x)))^ke^{-na(t(x))}\right]|\Omega(1,2)|\\
		&\leq & \Theta(1) \frac n{k!}(\frac{4}{3}n\pi r^3)^k e^{- n\pi r^3}|\Omega(1,2)|\\
		&=& \Theta(1)\frac{n(\log n+o(\log n))^k}{n+ (\log n)^{\Theta (1)}} \left(\frac{\log n}{n}\right)^{\frac{1}{3}}=o(1).\\
	\end{eqnarray*}
\end{proof}

\begin{figure}[htbp]
	\begin{center}
		\scalebox{0.45}[0.45]{
			\includegraphics{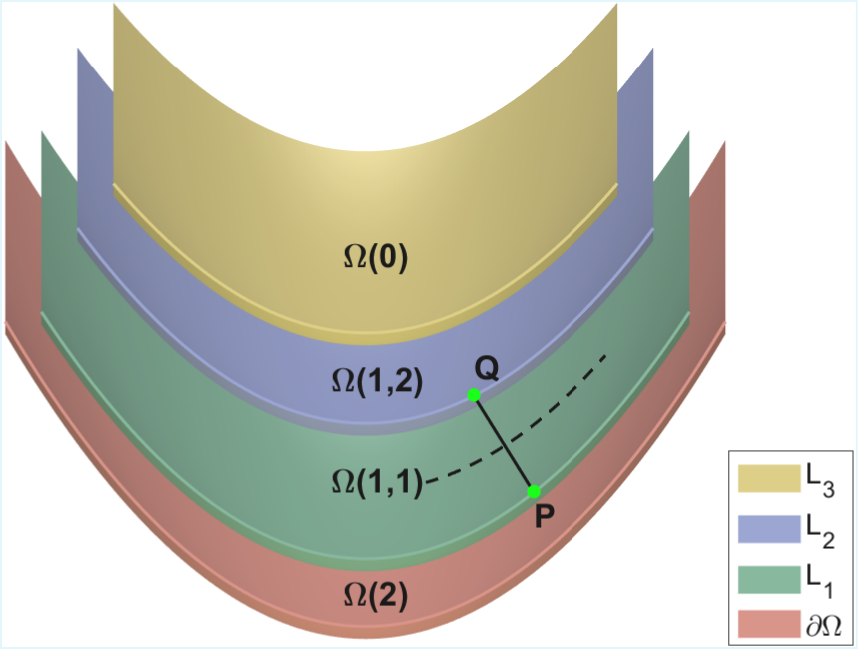}       }
		\caption{$\Omega$ is divided into four parts}\label{fig:Omega(1)}
	\end{center}
\end{figure}

The left work is to determine $n\int_{\Omega(1,1)} \psi^k_{n,r}(x)\mathrm{d}x$.
By  (\ref{eq:estimation-psi(n)}),
\begin{equation}\label{eq:Omega(1,1)-approximation}
	\begin{aligned}
		&\frac n{k!}\int_{\Omega(1,1)}\left(n|B(x,r_n)\cap\Omega|\right)^k e^{-n|B(x,r_n)\cap\Omega|}\mathrm{d}x\\
		\sim & \frac
		n{k!}\int_{\Omega(1,1)}(na(t(x)))^ke^{-na(t(x))}\mathrm{d}x.
	\end{aligned}
\end{equation}
We define
\begin{equation*}\label{eq:L1}
L_1=\{x\in\Omega\mid \mathrm{dist}(x,\partial \Omega)=(G(\partial\Omega)+1)r^2\},
\end{equation*}
\begin{equation*}\label{eq:L2}
L_2=\left\{x\in\Omega(1)\mid t(x)=\frac{r}{2}\right\},
\end{equation*}
\begin{equation*}\label{eq:L3}
L_3=\{x\in\Omega\mid \mathrm{dist}(x,\partial \Omega)=r\}.
\end{equation*}
Clearly, the subregion $\Omega(1,1)$ has boundaries
$L_1$ and $L_2$,   illustrated by
Figure~\ref{fig:Omega(1)}. If $x\in L_1$, then $t(x)\geq 0$ by a similar argument   in   previous discussions. Combing the fact $t(x)\leq (G(\partial\Omega)+1)r^2$, we know
 $$
 0\leq t(x)\leq (G(\partial\Omega)+1)r^2, \quad  x\in L_1.
$$
As $r$
tends to zero, both surfaces $L_1$ and $L_2$ approximate the
boundary $\partial \Omega$. So there exists $\epsilon_r$ such
that  the area of $L_1$ and $L_2$,
denoted by $\mathrm{Area}(L_1)$ and $\mathrm{Area}(L_2)$ respectively, satisfies:
$$
\mathrm{Area}(\partial\Omega)-\epsilon_r\leq \mathrm{Area}(L_1),\mathrm{Area}(L_2)\leq \mathrm{Area}(\partial\Omega),
$$
where $\epsilon_r\rightarrow 0$ as $r\rightarrow 0$. In addition, it
is obvious that $t(x)$   increases along the directed line segment
started from a point in $L_1$ to a point in $L_2$. See the line
segment $PQ$ in Figure~\ref{fig:Omega(1)}.

The integration on   $\Omega(1,1)$  can be bounded as
follows:
\begin{equation}\label{eq:Bounds-Integration-Omega(1,1)}
\begin{split}
&(\mathrm{Area}(\partial\Omega)-\epsilon_r)\frac n{k!}\int_{(G(\partial\Omega)+1)r^2}^{\frac r2}(na(t))^ke^{-na(t)}\mathrm{d}t\\
&\leq \frac n{k!}
\int_{\Omega(1,1)}(na(t(x)))^ke^{-na(t(x))}\mathrm{d}x\\
&\leq \mathrm{Area}(\partial\Omega)\frac n{k!}\int_0^{\frac r2}(na(t))^ke^{-na(t)}\mathrm{d}t.
\end{split}
\end{equation}
Inequalities (\ref{eq:Bounds-Integration-Omega(1,1)}) can be similarly proved, by an analogous  argument to that of  inequalities~(18) presented in~\cite{Ding2025-IEEETON}.

According to Lemma~\ref{lemma:a(t)}, we obtain that
\begin{eqnarray*}
&&\mathrm{Area}(\partial\Omega)\frac n{k!}\int_0^{\frac r2}(na(t))^ke^{-na(t)} \mathrm{d}t\\
&\sim&\mathrm{Area}(\partial\Omega)  \frac{e^{-\frac{2\xi}{3}}}{ \pi^{\frac{1}{3}}}\left(\frac{2}{3}\right)^k\frac{1}{k!}.
\end{eqnarray*}
Therefore,
\begin{equation*}
\begin{split}
0\leq &\epsilon_r\int_{(G+1)r^2}^{\frac{r}{2}}\frac
n{k!}(na(t))^ke^{-na(t)}\mathrm{d}t\\
\leq& \epsilon_r\int_{0}^{\frac{r}{2}}\frac
n{k!}(na(t))^ke^{-na(t)}\mathrm{d}t=o(1).
\end{split}
\end{equation*}
Furthermore, noticing that $2\pi r^3/3\leq a(t)\leq 2(\pi+\delta)r^3/3$
for any $t\in [0,(G(\partial\Omega)+1)r^2]$, where $\delta>0$ is a small real number. By simple calculation,
\begin{equation*}
\begin{aligned}
&\mathrm{Area}(\partial\Omega)\int_0^{(G+1)r^2}\frac n{k!}(na(t))^ke^{-na(t)}\mathrm{d}t =o(1).
\end{aligned}
\end{equation*}
So,
\begin{equation*}
\begin{split}
&( \mathrm{Area}(\partial\Omega)-\epsilon_r)\frac n{k!}\int_{(G(\partial\Omega)+1)r^2}^{\frac r2}(na(t))^ke^{-na(t)}\mathrm{d}t\\
=&\frac{n}{k!} \mathrm{Area}(\partial\Omega)\int_0^{\frac{r}{2}}(na(t))^ke^{-na(t)}\mathrm{d}t\\
&- \frac{n}{k!}\mathrm{Area}(\partial\Omega)\int_0^{(G(\partial\Omega)+1)r^2} (na(t))^ke^{-na(t)}\mathrm{d}t\\
&-\frac{n}{k!}\epsilon_r\int_{(G(\partial\Omega)+1)r^2}^{\frac{r}{2}} (na(t))^ke^{-na(t)}\mathrm{d}t\\
=& \mathrm{Area}(\partial\Omega)\frac {n}{k!}\int_0^{\frac{r}{2}}(na(t))^ke^{-na(t)}\mathrm{d}t+o(1).
\end{split}
\end{equation*}
Consequently, inequalities~(\ref{eq:Bounds-Integration-Omega(1,1)}) can be expressed as:
\begin{equation}\label{eq:Omega(1,1)-Inequality}
\begin{split}
& \mathrm{Area}(\partial\Omega)\frac {n}{k!}\int_0^{\frac{r}{2}}(na(t))^ke^{-na(t)}\mathrm{d}t+o(1)\\
=&( \mathrm{Area}(\partial\Omega)-\epsilon_r)\frac n{k!}\int_{(G(\partial\Omega)+1)r^2}^{\frac r2}(na(t))^ke^{-na(t)}\mathrm{d}t\\
\leq&\frac n{k!}\int_{\Omega(1,1)}(na(t(x)))^ke^{-na(t(x))}\mathrm{d}x\\
\leq&  \mathrm{Area}(\partial\Omega)\frac n{k!}\int_0^{\frac r2}(na(t))^ke^{-na(t)}\mathrm{d}t,
\end{split}
\end{equation}
which implies that
\begin{eqnarray*}
&&\frac n{k!}\int_{\Omega(1,1)}(na(t(x)))^ke^{-na(t(x))}\mathrm{d}x \\
&\sim &  \mathrm{Area}(\partial\Omega)\frac n{k!}\int_0^{\frac r2}(na(t))^ke^{-na(t)} \mathrm{d}t.
\end{eqnarray*}
Hence, by (\ref{eq:Omega(1,1)-approximation}) and Lemma~\ref{lemma:a(t)}, we have
\begin{eqnarray*}
&&\frac n{k!}\int_{\Omega(1,1)}\left(n|B(x,r_n)\cap\Omega|\right)^k e^{-n|B(x,r_n)\cap\Omega|}\mathrm{d}x\\
&\sim &\frac{n}{k!}\int_{\Omega(1,1)}(na(t(x)))^ke^{-na(t(x))}\mathrm{d}x\\
&\sim&  \mathrm{Area}(\partial\Omega)\frac {n}{k!}\int_0^{\frac r2}(na(t))^ke^{-na(t)}\mathrm{d}t\\
&\sim&  \mathrm{Area}(\partial\Omega)  \frac{e^{-\frac{2\xi}{3}}}{ \pi^{\frac{1}{3}}}\left(\frac{2}{3}\right)^k\frac{1}{k!}.
\end{eqnarray*}
Therefore, we have proved the following conclusion:
\begin{claim}\label{pro:Omega(1,1)}
$$
n\int_{\Omega(1,1)} \psi^k_{n,r}(x)\mathrm{d}x\sim
\mathrm{Area}(\partial\Omega)  \frac{e^{-\frac{2\xi}{3}}}{ \pi^{\frac{1}{3}}}\left(\frac{2}{3}\right)^k\frac{1}{k!}.
$$
\end{claim}
The four claims prove
\begin{eqnarray*}
&&\int_{\Omega}n\psi^k_{n,r}(x)\mathrm{d}x\\
&=&\left\{\int_{\Omega(0)}+\int_{\Omega(2)}+\int_{\Omega(1,2)}+\int_{\Omega(1,1)}\right\}
n\psi^k_{n,r}(x)\mathrm{d}x \\
&\sim&  \mathrm{Area}(\partial\Omega)  \frac{e^{-\frac{2\xi}{3}}}{ \pi^{\frac{1}{3}}}\left(\frac{2}{3}\right)^k\frac{1}{k!} = e^{-c}.
\end{eqnarray*}

\section{Proof of Proposition~\ref{pro:Explicit-Form} under condition $C_{\mathrm{II}}$}\label{sec:Pro-1}

For convenience and without the loss of generality,  we assume that boundary $\partial\Omega$ are smooth everywhere except at a curve $\sigma\subset\partial\Omega$  of length $\mathrm{Len}(\sigma)$. That is, the singularity part of the boundary is $\sigma$ and the smooth part is $\partial\Omega\backslash \sigma$.
We   use notation $G(\Omega\backslash\sigma)$ to represent the maximum curvature function of surface $\partial\Omega\backslash\sigma$, which can be similarly derived under the circumstance of  surface $\partial\Omega\backslash\sigma$, as Corollary~\ref{corollary:Bounds-B(x,r)-Omega-sigma} illustrates.

Recall  $\mathrm{Tube}_{2r}(\sigma), \mathrm{Tube}_{3r}(\sigma)$, the tubes containing curve $\sigma$  defined by (\ref{eq:tubes-sigma}), it is easy to see that
$$
| \mathrm{Tube}_{3r}(\sigma)|\leq (\mathrm{Len}(\sigma)+6r) \pi (3r)^2=\Theta(r^2).
$$
As a result, by Remark~\ref{rem:remark3} (iii),
\begin{equation}\label{eq:Tube(sigma)}
	\begin{aligned}
		& n\int_{ \mathrm{Tube}_{3r}(\sigma)} \psi^k_{n,r}(x)\mathrm{d}x \\
		\leq &\frac n{k!}\left(\frac{4}{3}n\pi r^3\right)^ke^{-\frac{4C_*}{3}n\pi r^3}  | \mathrm{Tube}_{3r}(\sigma)|\\
		\leq & \frac{\Theta(1)}{r}\left(n\pi r^3\right)^{k+1}e^{-\frac{4C_*}{3}n\pi r^3}=o(1),
	\end{aligned}
\end{equation}
where constant $C_*>\frac{1}{4}$ due to Condition~$C_{\mathrm{II}}$.

Let $  \Omega'= \Omega\backslash  \mathrm{Tube}_{3r}(\sigma)$, and  we
similarly define
$$
\Omega'(0)=\left\{x\in\Omega':\mathrm{dist}(x,\partial\Omega)\geq r\right\}, $$
$$
\Omega'(2)=\left\{x\in\Omega':\mathrm{dist}(x,\partial\Omega)\leq (G(\Omega\backslash\sigma)+1)r^2\right\},
$$
$$
\Omega'(1)=\Omega'\backslash (\Omega'(0)\cup\Omega'(2)).
$$
Therefore
$
\Omega=\Omega' \cup  \mathrm{Tube}_{3r}(\sigma)=\cup_{i=0}^2 \Omega'(i) \cup  \mathrm{Tube}_{3r}(\sigma).
$

\par Follow the proofs of Claim~\ref{pro:Omega(0)} and~\ref{pro:Omega(2)}, it is easy to verify
\begin{equation}\label{eq:Omega'(0)}
	\begin{aligned}
		n\int_{\Omega'(0) } \psi^k_{n,r}(x)\mathrm{d}x =n\int_{\Omega'(2) } \psi^k_{n,r}(x)\mathrm{d}x =o(1).
	\end{aligned}
\end{equation}
The following deals with $n\int_{\Omega'(1) } \psi^k_{n,r}(x)\mathrm{d}x$.
For any point $x\in \Omega'(1)$,
$$
(G(\partial\Omega\backslash\sigma)+1)r^2< \mathrm{dist}(x,\partial \Omega)<r<3r<\mathrm{dist}(x,\sigma).
$$
As discussed in Section~\ref{sec:B(x,r)-Omega-bound-sigma} and illustrated by Corollarly~\ref{corollary:Bounds-B(x,r)-Omega-sigma},
$
|B(x,r)\cap\Omega|\subset \Omega\backslash \mathrm{Tube}_{2r}(\sigma),
$
and
$$
a(t(x,r))\leq  	|B(x,r)\cap\Omega|\leq   a(t(x,r))+\Theta(r^4).
$$
Therefore, $\forall x\in \Omega'(1)$,
\begin{equation}\label{eq:local-1-1}
	\begin{split}
		& (1+o(1))(na(t(x)))^ke^{-na(t(x))}\\
	 \leq &(n|B(x,r)\cap \Omega|)^ke^{-n|B(x,r)\cap \Omega|}\\
		\leq& (1+o(1))^k(na(t(x)))^ke^{-na(t(x))}.
	\end{split}
\end{equation}
Furthermore, we can similarly set
$$
\Omega'(1,1)=\left\{x\in \Omega'(1):t(x)\leq \frac r2\right\},
$$
$$
\Omega'(1,2)=\Omega'(1)\setminus\Omega'(1,1),
$$
where $t(x)$ is similarly  determained according to Corollarly~\ref{corollary:Bounds-B(x,r)-Omega-sigma} as well as Lemma~\ref{lemma:Bounds-B(x,r)-Omega}.

By the same argument of the proof of Claim~\ref{pro:Omega(1,2)}, we will see
\begin{equation}\label{eq:Omega'(1,2)}
	n\int_{\Omega'(1,2)} \psi^k_{n,r}(x)\mathrm{d}x=o(1).
\end{equation}
In fact,
$|\Omega'(1,2)|\leq \mathrm{Area}(\partial\Omega)r\sim \mathrm{Area}(\partial\Omega) \left(\frac{\log n}{n}\right)^{\frac{1}{3}},
$
and  $
a(t(x))\geq \pi r^3
$, which is due to $t(x)\geq \frac{r}{2}$ for  $x\in \Omega'(1,2)$.
By inequalities~(\ref{eq:local-1-1}),
\begin{eqnarray*}
	&&n\int_{\Omega'(1,2)} \psi^k_{n,r}(x)\mathrm{d}x\\
	&=&\frac n{k!}\int_{\Omega'(1,2)}\left(n|B(x,r_n)\cap\Omega|\right)^k e^{-n|B(x,r_n)\cap\Omega|}\mathrm{d}x\\
	&\leq&\frac n{k!}\left[(1+o(1))^k(na(t(x)))^ke^{-na(t(x))}\right]|\Omega'(1,2)|\\
	&\leq & \Theta(1) \frac n{k!}\left(\frac{4}{3}n\pi r^3\right)^k e^{- n\pi r^3}|\Omega'(1,2)|\\
	&=& \Theta(1)\frac{n(\log n+o(\log n))^k}{n+ (\log n)^{\Theta (1)}} \left(\frac{\log n}{n}\right)^{\frac{1}{3}}=o(1).
\end{eqnarray*}

Notice that $\mathrm{Area}( \mathrm{Tube}_{3r}(\sigma))\leq \Theta(1)r$, and
\begin{equation*}
	\begin{split}
		\mathrm{Area}(\partial\Omega\backslash \mathrm{Tube}_{3r}(\sigma))=\mathrm{Area}(\partial\Omega)+o(1).
	\end{split}
\end{equation*}
Then similarly to (\ref{eq:Omega(1,1)-Inequality}), we have
\begin{equation*}\label{eq:Omega(1,1)-Inequality-2}
	\begin{split}
		& \mathrm{Area}(\partial\Omega\backslash \mathrm{Tube}_{3r}(\sigma))\frac {n}{k!}\int_0^{\frac{r}{2}}(na(t))^ke^{-na(t)}\mathrm{d}t+o(1)\\
		\leq&\frac n{k!}\int_{\Omega'(1,1)}(na(t(x)))^ke^{-na(t(x))}\mathrm{d}x\\
		\leq&  \mathrm{Area}(\partial\Omega\backslash\mathrm{Tube}_{3r}(\sigma))\frac n{k!}\int_0^{\frac r2}(na(t))^ke^{-na(t)}\mathrm{d}t+o(1).
	\end{split}
\end{equation*}
This further implies that
\begin{equation}\label{eq:Omega'(1,1)}
	\begin{split}
		n\int_{\Omega'(1,1)} \psi^k_{n,r}(x)\mathrm{d}x\sim&  \mathrm{Area}(\partial\Omega)\frac {n}{k!}\int_0^{\frac r2}(na(t))^ke^{-na(t)}\mathrm{d}t\\
		\sim&  \mathrm{Area}(\partial\Omega)  \frac{e^{-\frac{2\xi}{3}}}{ \pi^{\frac{1}{3}}}\left(\frac{2}{3}\right)^k\frac{1}{k!}.
	\end{split}
\end{equation}
Combining (\ref{eq:Tube(sigma)}), (\ref{eq:Omega'(0)}), (\ref{eq:Omega'(1,2)}) and (\ref{eq:Omega'(1,1)}), we know
\begin{equation*}
	\begin{split}
		&\int_{\Omega}n\psi^k_{n,r}(x)\mathrm{d}x\\
		=&\int_{\mathrm{Tube}_{3r}(\sigma)}	n\psi^k_{n,r}(x)\mathrm{d}x \\		
		&+\left\{\int_{\Omega'(0)}+\int_{\Omega'(2)}+\int_{\Omega'(1,2)}+\int_{\Omega'(1,1)}\right\}
		n\psi^k_{n,r}(x)\mathrm{d}x \\
		\sim&  \mathrm{Area}(\partial\Omega)  \frac{e^{-\frac{2\xi}{3}}}{ \pi^{\frac{1}{3}}}\left(\frac{2}{3}\right)^k\frac{1}{k!} = e^{-c}.
	\end{split}
\end{equation*}

\begin{remark}\label{remark-2:Condition3}
 The techniques dealing with 	Condition~$C_\mathrm{II}$ can also applied to  	Condition~$C_\mathrm{III}$. The finite number of singurality points, consists of singurality set  $\sigma$,  can be covered by a finite number of balls centered at these points with the same radius $3r$. Let  $\mathrm{Tube}_{3r}(\sigma)$ be the union of these balls, then its volume
  $
 | \mathrm{Tube}_{3r}(\sigma)|\leq  \Theta(r^3).
 $
The boundary requirement   $\liminf_{r\rightarrow 0} \inf_{x\in\Omega}\frac{\mathrm{Vol}(B(x,r)\cap \Omega )}{\mathrm{Vol}(B(x,r))}>\frac{1}{4}$
can be correspondingly relaxed to  $\liminf_{r\rightarrow 0} \inf_{x\in\Omega}\frac{\mathrm{Vol}(B(x,r)\cap \Omega )}{\mathrm{Vol}(B(x,r))}>0$. The other proof precedure for  Condition~$C_\mathrm{III}$ is the same to   Condition~$C_\mathrm{II}$.

\end{remark}

\section{Proof of Proposition~\ref{pro:conclusion-1} under condition~$C_{\mathrm{I}}$}\label{sec:Pro-2}
To prove  Proposition~\ref{pro:conclusion-1}, we follow the way,   instead to prove  its Poissonized version.
Then by the standard de-Poissonized technique to see that Proposition~\ref{pro:conclusion-1} holds.

\subsection{Poissonized version of Proposition~\ref{pro:conclusion-1}}

 A homogeneous Poisson point process $\mathcal{P}_n$ distributed over region $\Omega\subset \mathbb{R}^d\ (d\geq 1)$ with intensity $n|\Omega|=n$, satisfies that
	
	 \begin{enumerate}
	 	\item \textbf{Poisson Distribution}:   	 	For any  measurable set $A \subset \Omega$, the number of points in $A$,  $\mathcal{P}_n(A) = \#(\mathcal{P}_n\cap A)$, follows a Poisson distribution  with intensity $n|A|$:
	 	$$
	 	\Pr\left\{\mathcal{P}_n(A)=j\right\}=\frac{e^{-n|A|}(n|A|)^j}{j!},\quad j\geq0.
	 	$$
	 	\item \textbf{Independence}:
	 	For any finite collection of disjoint measuralbe sets $A_1, \dots, A_K$, the random variables $\mathcal{P}_n(A_1), \dots, \mathcal{P}_n(A_K)$ are mutually independent.
	 \end{enumerate}
	 In addition, given the condition that there are exactly $j$ points in $A$,   these $j$ points are independently and uniformly distributed in $A$.		
	These properties make a uniform $n$-point process  $\chi_n$ can be well approximated by a homogeneous Poisson point process $\mathcal{P}_n$.

   Moreover, let $\{X_i\}_{i=1}^\infty$ be a sequence of independent random $d$-dimensional vectors with common probability distribution $F$, whose probability density function $f$ is continuous almost everywhere. Define the Poisson point process
$
\mathcal{P}_\lambda = \{X_1, X_2, \dots, X_{N_\lambda}\},
$
where $N_\lambda$ is a Poisson random variable, independent of the sequence $\{X_i\}_{i=1}^\infty$.
Then \( G(\mathcal{P}_\lambda, r) \) is a geometric graph on the point set \( \mathcal{P}_\lambda \) with connection radius \( r \); that is, any two points whose Euclidean distance is at most \( r \) are connected by an edge.     By $\mathrm{Po}(\beta)$ we denote a Poisson random variable with parameter $\beta$.

\begin{lemma}\label{lem:Poisson-approximation}(Theorem~6.7, \cite{Penrose-RGG-book})
	Let $W'_{k,\lambda}(r)$ be the number of vertices of degree $k$ in $G(\mathcal{P}_\lambda, r)$.
	For $x\in \mathbb{R}^d$, denote $B_x=B(x,r)$ and $\mathcal{P}_\lambda^x=\mathcal{P}_\lambda\cup\{x\}$. Define
	\begin{eqnarray*}
		J_1  &=&\lambda^2\int_{R^d}\mathrm{d}F(x)\\
		&\times&\int_{B(x,3r)}\mathrm{d}F(y)\Pr\{\mathcal{P}_\lambda(B_x)=k\}\Pr\{\mathcal{P}_\lambda(B_y)=k\},
	\end{eqnarray*}
	\begin{eqnarray*}
		J_2&=&\lambda^2\int_{R^d}\mathrm{d}F(x)\\
		&\times&\int_{B(x,3r)}\mathrm{d}F(y)\Pr\{ \{\mathcal{P}^y_\lambda(B_x)=k\}\cap \{\mathcal{P}^x_\lambda(B_y)=k\}\}.
	\end{eqnarray*}
	Then
	$$
	d_{\mathrm{TV}}(W'_{k,\lambda}(r), \mathrm{Po}(E(W'_{k,\lambda}(r))))\leq 3 (J_1+J_2),
	$$
	where $d_{\mathrm{TV}}$ refers to total variation distance.
\end{lemma}
This lemma indicates that \( W'_{k,\lambda}(r) \) can be approximated in distribution by a Poisson random variable with parameter \( \mathbb{E}[W'_{k,\lambda}(r)] \), provided that \( J_1 + J_2 \to 0 \).

For the uniform \( n \)-point process \( \chi_n \) over \( \Omega \) considered in this article---where the density is uniform, i.e., \( f(x) = 1 \) for \( x \in \Omega \)---Palm theory (Theorem~1.6, \cite{Penrose-RGG-book}) yields
\[
\mathbb{E}[W'_{k,n}(r)] = n \int_{\Omega} \psi^k_{n,r}(x) f(x)\, \mathrm{d}x = n \int_{\Omega} \psi^k_{n,r}(x)\, \mathrm{d}x \sim e^{-c}.
\]
The asymptotic relation \( n \int_{\Omega} \psi^k_{n,r}(x)\, \mathrm{d}x \sim e^{-c} \) is established in Proposition~\ref{pro:Explicit-Form}.

%

That is, the random variable $ W'_{k,\lambda}(r) $, which counts the number of vertices of degree $ k $, is asymptotically Poisson distributed with parameter $ e^{-c} $, i.e.,
$ W'_{k,\lambda}(r) \xrightarrow{d} \mathrm{Po}(e^{-c}) $,
as guaranteed by Lemma~\ref{lem:Poisson-approximation} (provided that $ J_1 + J_2 \to 0 $).

Consequently, the proof of Proposition~\ref{pro:conclusion-1},
\[
\lim_{n \to \infty} \Pr\!\left\{ \rho(\chi_n; \delta \geq k+1) \leq r_n \right\} = e^{-e^{-c}},
\]
can be carried out in two steps.
First, we establish the following Poissonized counterpart, Proposition~\ref{pro:Poisson-Version}, which corresponds to Proposition~\ref{pro:conclusion-1} under a Poisson point process.
Second, we apply a standard de-Poissonization argument to transfer the result back to the original   setting. Since this de-Poissonization technique is classical and well-documented, we omit its details here and refer the reader to~\cite{Penrose-k-connectivity} for a thorough treatment.

\begin{proposition}\label{pro:Poisson-Version}
	Under the assumptions of Theorem~\ref{thm:Main3D},
	\begin{equation}\label{eq:concusion-1-Poisson}
		\lim_{n\rightarrow\infty}\Pr\left\{\rho(\mathcal{P}_n;\delta\geq k+1)\leq r_n\right\} =\exp\left(-e^{-c}\right),
	\end{equation}
	where $\mathcal{P}_n$ is a homogeneous Poisson point process of intensity $n$ (i.e., $n|\Omega|$) distributed over unit-volume  region $\Omega$.
\end{proposition}

In the following, we establish the Poissonized version~\eqref{eq:concusion-1-Poisson} for a three-dimensional region $\Omega$ satisfying Condition~$C_{\mathrm{I}}$.
Before presenting the proof, we first provide a general estimate for Poisson point processes in the next subsection.

\subsection{A conclusion on Poisson process}
Denote
$$
v_{y\backslash x}=|B(y,r)\cap \Omega|-|B(x,r)\cap B(y,r)\cap \Omega|.
$$

\begin{proposition}\label{proposition1}Under the assumptions of Condition~$C_{\mathrm{I}}$ in Theorem~\ref{thm:Main3D}, let $x \in \Omega$ and let $Z_3$ be a Poisson random variable with mean $n v_{y \setminus x}$ as defined above. Then for $ j\geq 0$,
	\begin{equation}\label{eq:Pro-5}
		\frac{1}{n\pi r^3}\int_{\Omega\cap
			B(x,r)}n\Pr(Z_3=k-1-j)\mathrm{d}y=o(1).
	\end{equation}
\end{proposition}

\begin{proof} Here we  only  prove the case of $j=0$.  For $j>0$, the
	proof is similar and omitted.
	Let $
	d_0=\left(\frac{4r}{n^{\frac23}\pi^{\frac23}}\right)^{\frac13}$, then for any sufficiently large $n$,
	$$
	0<d_0= \Theta\left(\frac{(\log n)^\frac{1}{9}}{n^{\frac13}}\right)<r,
	$$
	and
	$
	\frac{n \pi d_0^3}{4}=(n\pi r^3)^{\frac13}.
	$
	$\forall x\in \Omega$, let
	$$
	\Gamma_1(x)=\left\{y\in B(x,r)\cap \Omega: \mathrm{dist}(y,x)\leq d_0\right\},
	$$
		\begin{equation}
		\begin{split}
			&\Gamma_2(x)=(B(x,r)\cap\Omega)\bigcap\\
			&\left\{y:   \mathrm{dist}(y,x)\geq d_0, \mathrm{dist}(y,\partial \Omega)> (G(\partial\Omega)+1)r^2\right\},
		\end{split}
	\end{equation}

	$$
	\Gamma_3(x)=\left\{y\in B(x,r)\cap \Omega: \mathrm{dist}(y,\partial \Omega)\leq (G(\partial\Omega)+1)r^2\right\}.
	$$
 Here  $G(\partial\Omega)$ 	is the curvature constant given in Lemma~\ref{lemma:G(Omega)}.  Obviously,
	$
	B(x,r)\cap\Omega\subset\Gamma_1(x)\cup\Gamma_2(x)\cup\Gamma_3(x),
	$
	as Figure~\ref{fig:Ball-in-Gammas} illustrates.


\begin{figure}[htbp]
	\begin{center}
		\scalebox{0.45}[0.45]{
			\includegraphics{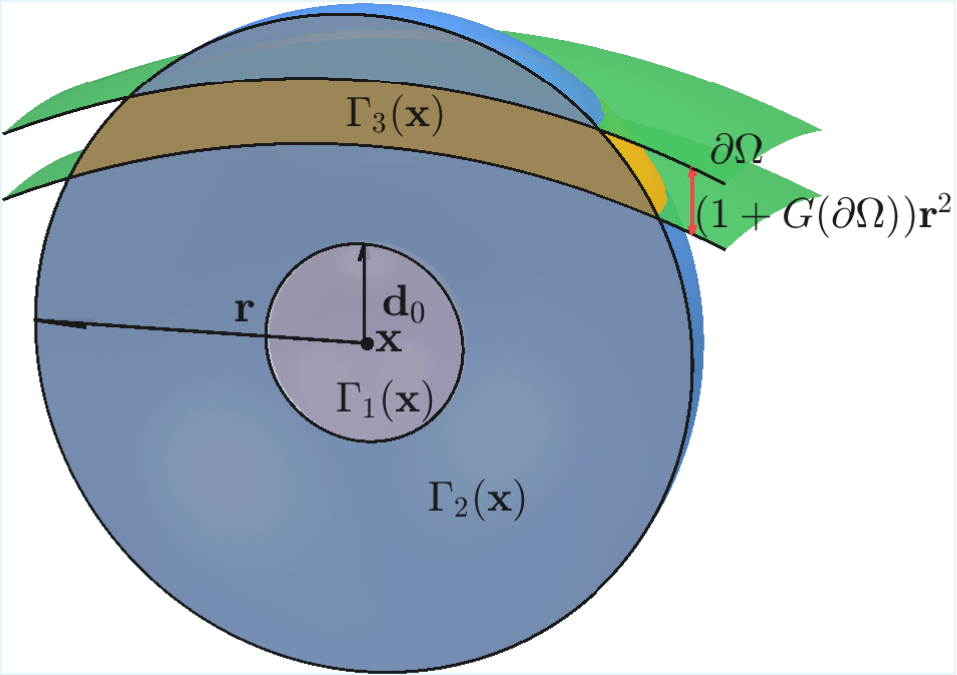}}
		\caption{ $B(x,r)\cap\Omega\subset\Gamma_1(x)\cup\Gamma_2(x)\cup\Gamma_3(x)$}\label{fig:Ball-in-Gammas}
	\end{center}
\end{figure}

\begin{figure}[htbp]
	\begin{center}
		\scalebox{0.43}[0.43]{
			\includegraphics{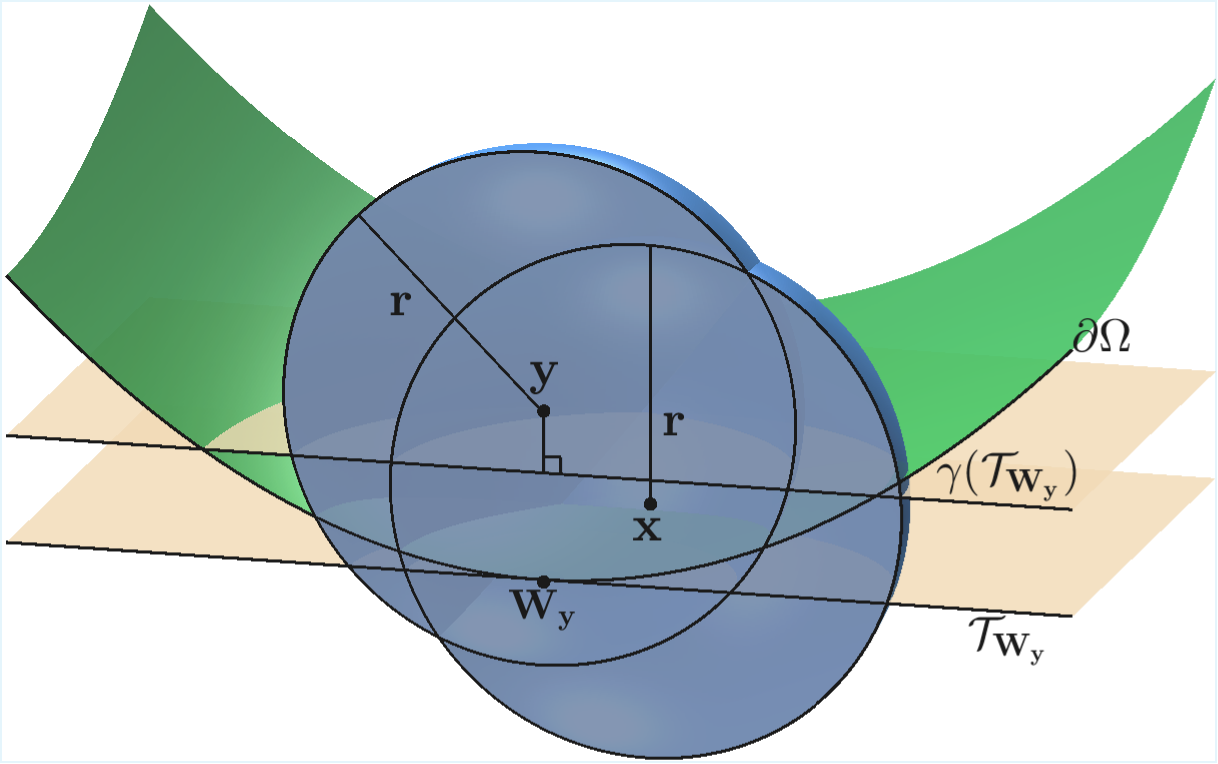}}
		\caption{Illustration of $A_2=o(1)$}\label{fig:lemma2-Application}
	\end{center}
\end{figure}

For $i=1,2,3$, let
\begin{equation}\label{eq:A123}
A_i=\frac{1}{n\pi r^3} \int_{\Gamma_i(x)}n\Pr\{Z_3=k-1\}\mathrm{d}y.
\end{equation}
We will prove $A_i=o(1), i=1,2,3$, in the following three steps.

\par Step 1: to prove $A_1=o(1)$.
\begin{eqnarray*}
A_1&=&\frac{1}{n\pi r^3} \int_{\Gamma_1(x)}n\Pr\{Z_3=k-1\}\mathrm{d}y\\
&\leq& \frac{1}{n\pi r^3} n |\Gamma_1(x)|
 \leq \Theta(1)\frac{ d_0^3}{ r^3}  =o(1).
\end{eqnarray*}

\par Step 2: to prove $A_2=o(1)$. For any $y\in \Gamma_2(x)$,
$\mathrm{dist}(y,\partial \Omega)\geq(G(\partial\Omega)+1)r^2$.
Let $W_y$ be the nearest point on $\partial\Omega$ to $y$, and   $\mathcal{T}_{W_y}$ be the tangent plane of $\partial\Omega$ at point $W_y$.
Then by the discussion in Section~\ref{sec:B(x,r)-Omega-Condition-I} and  (\ref{eq:B(x,r)-Omega-subset}), $B(y,r)\cap\Omega$ contains the major    segment of sphere $B(y,r)$  cut by a  plane $\gamma(\mathcal{T}_{W_y}) \pll \mathcal{T}_{W_y}$, i.e.,
\begin{equation} \label{eq:segment-ball-contain}
	B(y,r)|_{\gamma(\mathcal{T}_{W_y})}\subset B(y,r)\cap\Omega.
\end{equation}
See Figure~\ref{fig:lemma2-Application}.
This means that   more than one half of  $B(y,r)$ falls in $\Omega$.
 Consequently,
\begin{equation} \label{eq:vy-x}
	v_{y\backslash x}=|B(y,r)\cap \Omega|-|B(x,r)\cap B(y,r)\cap \Omega|\geq \frac{V^{\ast}(d_0)}{2},
\end{equation}	
where  $V^{\ast}(d_0)=\frac{\pi d_0^3}{4}$ is given by Lemma~\ref{lem:shadow-low-bound}. Therefore,
$$
nv_{y\backslash x}\geq\frac{1}{2}nV^{\ast}(d_0)=\frac{n}{2}\frac{\pi d_0^3}{4}=\frac{1}{2}(n\pi r^3)^{\frac13}.
$$

\begin{eqnarray*}
A_2
&=&\frac{1}{n\pi r^3} \int_{\Gamma_2(x)}n\Pr\{Z_3=k-1\}\mathrm{d}y\\
&\leq&\frac{1}{n\pi r^3} \int_{\Gamma_2(x)}\frac{n(n\pi r^3)^{k-1}e^{-\frac 12\left(n\pi r^3\right)^{\frac13} }}{(k-1)!} \mathrm{d}y\\
&\leq&\frac{1}{n\pi r^3}
\frac{(n\pi r^3)^{k}e^{-\frac {1}{2}\left(n\pi r^3\right)^{\frac13} }}{(k-1)!}
 \sim \frac{o(1)}{n\pi r^3} =o(1).
\end{eqnarray*}

\par Step 3: to prove $A_3=o(1)$.
Notice that $\Gamma_3(x)$ falls in a region with the width less than $(G(\partial\Omega)+1)r^2$, we have
$     |\Gamma_3(x)|\leq  \mathrm{Area}(\partial B(x,r)) (G(\partial\Omega)+1)r^2=\Theta(1)r^4.  $
\begin{eqnarray*}
	A_3&=&\frac{1}{n\pi r^3} \int_{\Gamma_3(x)}n\Pr\{Z_3=k-1\}\mathrm{d}y\\
	&\leq & \frac{1}{\pi r^3} |\Gamma_3(x)|
	 \leq    \Theta(r)  =o(1).
\end{eqnarray*}
Finally, we have
\begin{eqnarray*}
	&&\frac{1}{n\pi r^3} \int_{\Omega\cap B(x,r)}n\Pr\{Z_3=k-1\}\mathrm{d}y\\
	&\leq& A_1+A_2+A_3\\
	&\leq &\Theta\left(\frac{d_0^3}{r^3}+\frac{1}{n\pi r^3}+r\right)
	 \leq  \Theta\left(\frac{1}{(\log n)^\frac{2}{3}}\right) =o(1).
\end{eqnarray*}

The proposition is therefore proved.

\end{proof}

\subsection{ Proof  of Proposition~\ref{pro:Poisson-Version} under   condition~$C_{\mathrm{I}}$ }\label{sec:sub-Proof-Poisson-Version}

We follow Penrose's approach and framework to prove Proposition~\ref{pro:Poisson-Version} under condition~$C_{\mathrm{I}}$. The proof of such results in the two-dimensional case can be found in~\cite{Ding2025-IEEETON}. Here we deal with the three-dimensional case.

For given $n$, $x,y\in \Omega\subset\mathbb{R}^3$, let
$$
v_x=|B(x,r)\cap \Omega |, \quad v_y=|B(y,r)\cap \Omega |, $$
$$ v_{x,y}=|B(x,r)\cap B(y,r)\cap \Omega|,$$
$$v_{x\backslash y}=v_x-v_{x,y},v_{y\backslash x}=v_y-v_{x,y}.
$$
Define $I_i=I_i(n)(i=1,2,3)$ as follows:
$$
I_1=n^2\int_{\Omega} \mathrm{d}x\int_{\Omega\cap B(x,3r)}\mathrm{d}y \psi^k_{n,r}(y)\psi^{k}_{n,r}(x),
$$
$$
I_2=n^2\int_{\Omega} \mathrm{d}x\int_{\Omega\cap B(x, r)}\mathrm{d}y \Pr\{Z_1+Z_2=Z_1+Z_3=k-1\},
$$
$$
I_3=n^2\int_{\Omega} \mathrm{d}x\int_{\Omega\cap B(x,3r)\backslash B(x,r)}\mathrm{d}y \Pr\{Z_1+Z_2=Z_1+Z_3=k\},
$$
where $Z_1,Z_2,Z_3$ are independent Poisson variables with means
$nv_{x,y},nv_{x\backslash y},nv_{y\backslash x}$ respectively. As Penrose has pointed out in~\cite{Penrose-k-connectivity},
by an argument similar to that of Section~7 of~\cite{Penrose97:longest-edge}, to prove Proposition~\ref{pro:Poisson-Version} it suffices to prove that
$I_1,I_2,I_3\rightarrow 0$ as $n\rightarrow \infty$.
These conclusions will be proved by the following three claims.

\begin{claim}\label{claim:I-1}
	$I_1\rightarrow 0$ as $n\rightarrow \infty$.
\end{claim}
\begin{proof}
	By the assumptions and Lemma~\ref{lemma:C2-Boundary-volume}, there exists a constant $C>0$  such that $C\pi r^3\leq |B(x,r)\cap\Omega|\leq \frac{4}{3}\pi r^3$.
	
%

	\begin{eqnarray*}
		&&I_1\\
		&=&n^2\int_{\Omega} \mathrm{d}x\int_{{\Omega}\cap B(x,3r)}\mathrm{d}y \psi^k_{n,r}(y)\psi^{k}_{n,r}(x)\\
		&=& \int_{\Omega} n\psi^{k}_{n,r}(x)\mathrm{d}x \int_{B(x,3r)\cap \Omega}n\psi^{k}_{n,r}(y)\mathrm{d}y\\
		&\leq& \Theta(1)\int_{\Omega} n\psi^{k}_{n,r}(x)\mathrm{d}x \int_{B(x,3r)\cap \Omega}\frac{1}{k!} n(n\pi r^3)^ke^{-Cn\pi r^3}\\
		&\leq & \Theta(1)\frac{n}{k!}|B(x,3r)|(n\pi r^3)^{k}e^{-Cn\pi r^3}\int_{\Omega} n\psi^{k}_{n,r}(x)\mathrm{d}x\\
		&=& \Theta(1)(n\pi r^3)^{k+1}e^{-Cn\pi r^3}\int_{\Omega} n\psi^{k}_{n,r}(x)\mathrm{d}x\\
		&\sim& \Theta(1)(n\pi r^3)^{k+1}e^{-Cn\pi r^3}=o(1).
	\end{eqnarray*}
\end{proof}

\begin{claim}\label{claim:I-2}
	$I_2\rightarrow 0$, as $n\rightarrow \infty$.
\end{claim}
\begin{proof}
	It is obvious that
	$Z_1+Z_2$ and $Z_3$ are Poisson variables with means $nv_x$ and
	$nv_{y\backslash x}$ respectively.
	Notice that
	\begin{eqnarray*}
		&&\Pr\{Z_1+Z_3=k-1| Z_1+Z_2=k-1\}\\
		&\leq &\sum_{j=0}^{k-1}\Pr\{Z_3=k-1-j |Z_1+Z_2=k-1\}\\
		&=&\sum_{j=0}^{k-1}\Pr\{Z_3=k-1-j\}.
	\end{eqnarray*}
	
	For any $x\in\Omega$, $\psi^{k}_{n,r}(x)=\frac{n|B(x,r)\cap \Omega|}{k}\psi^{k-1}_{n,r}(x),$
	by Lemma~\ref{lemma:C2-Boundary-volume}, there exists a constant $C>0$ such that
	$C\pi r^3\leq |B(x,r)\cap \Omega|$, and thus
	$ \psi^{k-1}_{n,r}(x)\leq \frac{k}{Cn\pi r^3}\psi^{k}_{n,r}(x).$
	Therefore,
	\begin{equation}\label{eq:I-2}
    \begin{aligned}
		 &I_2\\
		 =&n^2\int_{\Omega} \mathrm{d}x\int_{\Omega\cap B(x,r)}\mathrm{d}y \Pr\{Z_1+Z_2=Z_1+Z_3=k-1\}\\
		 =&n^2\int_{\Omega} \mathrm{d}x\int_{\Omega\cap B(x,r)}\mathrm{d}y
		\Pr\{Z_1+Z_2=k-1\}\\
		 & \Pr\{Z_1+Z_3=k-1|Z_1+Z_2=k-1\}\\
		 =& n \int_{\Omega}\Pr\{Z_1+Z_2=k-1\} \mathrm{d}x\\
		 &\cdot\int_{\Omega\cap B(x,r)}n\Pr\{Z_1+Z_3=k-1|Z_1+Z_2=k-1\}\mathrm{d}y\\
		 \leq & n\int_{\Omega}\Pr\{Z_1+Z_2=k-1\} \mathrm{d}x\\
		 & \cdot\left\{\sum_{j=0}^{k-1}\int_{\Omega\cap B(x,r)}n\Pr\{Z_3=k-1-j\}\mathrm{d}y\right\}\\
		 =& \sum_{j=0}^{k-1} \int_{\Omega}n\psi^{k-1}_{n,r}(x) \mathrm{d}x \int_{\Omega\cap B(x,r)}n\Pr\{Z_3=k-1-j\}\mathrm{d}y\\
		 \leq& { \frac{k}{C}  \sum_{j=0}^{k-1}
			\frac{1}{n\pi r^3}\int_{\Omega}n\psi^{k}_{n,r}(x) \mathrm{d}x }\\
		 &\cdot\int_{\Omega\cap B(x,r)}n\Pr\{Z_3=k-1-j\}\mathrm{d}y\\
		 =&o(1).
    \end{aligned}
	\end{equation}
	The last equation holds due to Proposition~\ref{pro:Explicit-Form}  and~\ref{proposition1}.  	
\end{proof}

Similarly, we can prove
\begin{claim}\label{claim:I-3}
	$I_3\rightarrow 0$, as $n\rightarrow \infty$.
\end{claim}

Therefore, by these three claims,  we know the Poissonized version (\ref{eq:concusion-1-Poisson}) holds, which leads to
Proposition~\ref{pro:conclusion-1} by a de-Poissonized techinique as we have mentioned:
$$\lim_{n\rightarrow\infty}\Pr\left\{\rho(\chi_n;\delta\geq k+1)\leq r_n\right\} =e^{-e^{-c}}.$$

\section{Proof of Proposition~\ref{pro:conclusion-1} under condition~$C_\mathrm{II}$}\label{sec:Pro-2}
Under the condition~$C_\mathrm{II}$ in Theorem~\ref{thm:Main3D}, we will prove
$I_i\;(i=1,2,3)$ defined in Section~\ref{sec:sub-Proof-Poisson-Version} tend to zero. According to  condition~$C_\mathrm{II}$,
$$\frac{1}{3}\pi r^3 < \frac{4C_*}{3}\pi r^3\leq |B(x,r)\cap\Omega|\leq \frac{4}{3}\pi r^3
$$
 where $C_*>\frac{1}{4}$.
So Claim~\ref{claim:I-1} holds under condition~$C_\mathrm{II}$, i.e., $I_1\rightarrow 0, n\rightarrow\infty$. The following is to prove
$I_2\rightarrow 0, n\rightarrow\infty$. The convergence of $I_3$ can be similarly proved.

According to inequality~(\ref{eq:I-2}), it  sufficiently to prove
\begin{equation}\label{eq:I-2-main}
			\frac{1}{n\pi r^3}  \int_{\Omega}n\psi^{k}_{n,r}(x) \mathrm{d}x
		  \int_{\Omega\cap B(x,r)}n\Pr\{Z_3=k-1-j\}\mathrm{d}y=o(1).
\end{equation}
Here we  only  prove the case of $j=0$.
First, we need to prove that Proposition~\ref{proposition1} or (\ref{eq:Pro-5}) also holds
under   condition~$C_\mathrm{II}$.

\begin{eqnarray*}
	\frac{1}{n\pi r^3} \int_{\Omega\cap B(x,r)}n\Pr\{Z_3=k-1\}\mathrm{d}y\leq A_1+A_2+A_3  =o(1).
\end{eqnarray*}

\par By (\ref{eq:Tube(sigma)}), $\int_{\mathrm{Tube}_{3r}(\sigma)}n\psi^k_{n,r}(x)\mathrm{d}x=o(1)$, so
 \begin{equation}\label{eq:Pro-2-CII-a}
 \begin{aligned}
&\frac{1}{n\pi r^3}\int_{\mathrm{Tube}_{3r}(\sigma)}n\psi^k_{n,r}(x)\mathrm{d}x\int_{\Omega\cap
B(x,r)}n\Pr\{Z_3=k-1\}\mathrm{d}y\\
\leq & \frac{o(1)}{n\pi r^3} n|\Omega\cap B(x,r)| =o(1).
\end{aligned}
\end{equation}	

All left work is to prove
\begin{equation*}
\begin{aligned}
			&\frac{1}{n\pi r^3}  \int_{\Omega\backslash \mathrm{Tube}_{3r}(\sigma)}n\psi^{k}_{n,r}(x) \mathrm{d}x
		  \int_{\Omega\cap B(x,r)}n\Pr\{Z_3=k-1\}\mathrm{d}y\\
&=o(1).
\end{aligned}
\end{equation*}

In the following, we always assume $x\in\Omega\backslash\mathrm{Tube}_{3r}(\sigma)$.
Let
$$
	\Gamma'_1(x)=\left\{y\in B(x,r)\cap \Omega: \mathrm{dist}(y,x)\leq d_0\right\},
$$
\begin{equation*}
		\begin{split}
			&\Gamma'_2(x)=(B(x,r)\cap\Omega)\bigcap\\
			&\left\{ y:   \mathrm{dist}(x,y)\geq d_0, \mathrm{dist}(y,\partial \Omega)> (G(\partial\Omega\backslash\sigma)+1)r^2\right\},
		\end{split}
\end{equation*}
$$
	\Gamma'_3(x)=\left\{y\in B(x,r)\cap \Omega: \mathrm{dist}(y,\partial \Omega)\leq (G(\partial\Omega\backslash\sigma)+1)r^2\right\}.
$$
Here $d_0=\left(\frac{4r}{n^{\frac23}\pi^{\frac23}}\right)^{\frac13}$ as mentioned in the previous section.
Curvature constant $G(\partial\Omega\backslash \sigma)$ is   given by Corollary~\ref{corollary:Bounds-B(x,r)-Omega-sigma}.

\par Clearly,	$B(x,r)\cap\Omega\subset\Gamma'_1(x)\cup\Gamma'_2(x)\cup\Gamma'_3(x)$. Define
\begin{eqnarray*}
A'_i=\frac{1}{n\pi r^3} \int_{\Gamma'_i(x)}n\Pr\{Z_3=k-1\}\mathrm{d}y, \quad i=1,2,3.
\end{eqnarray*}
By similar arguments to the above section, we can easy prove  $A'_i=o(1),   i=1,3$.  The following is to determine
$A'_2=o(1)$.

\par
For any $y\in \Gamma'_2(x)$,
$\mathrm{dist}(y,\partial \Omega)\geq(G(\partial\Omega\backslash\sigma)+1)r^2$, and
 $\mathrm{dist}(y,\sigma)>2r$  due to  $y\in B(x,r)\cap \Omega$ and  $\mathrm{dist}(x,\sigma)>3r$.
 Let
$W_y$ be the nearest point on $\partial\Omega$ to $y$, then
\begin{equation*}
	\begin{aligned}
		\mathrm{dist}(W_y,\sigma)\geq & \mathrm{dist}(x,\sigma)-\mathrm{dist}(x,y) -\mathrm{dist}(y,W_y)\\
		> & 3r-r-(G(\partial\Omega\backslash\sigma)+1)r^2\\
		=& (2+o(1))r.
	\end{aligned}
\end{equation*}
That is to say, both $y$ and $W_y$ are suffciently far from the siguralirty part of the boundary.
 Then similarly to (\ref{eq:segment-ball-contain}),  as discussed in Section~\ref{sec:B(x,r)-Omega-Condition-I},  the major    segment of sphere $B(y,r)$ cut by a  plane $\gamma(\mathcal{T}_{W_y}) \pll \mathcal{T}_{W_y}$, is contained in  $B(y,r)\cap\Omega$, i.e. $ 	 B(y,r)|_{\gamma(\mathcal{T}_{W_y})}\subset B(y,r)\cap\Omega$.
So more than one half of  $B(y,r)$  falls in $\Omega$.
Again, similarly to (\ref{eq:vy-x}), we obtain
\begin{equation*}
\begin{aligned}
v_{y\backslash x}=&|B(y,r)\cap \Omega|-|B(x,r)\cap B(y,r)\cap \Omega|\geq \frac{V^{\ast}(d_0)}{2},
\end{aligned}
\end{equation*}
and
$
nv_{y\backslash x}\geq\frac{1}{2}nV^{\ast}(d_0)=\frac{n}{2}\frac{\pi d_0^3}{4}=\frac{1}{2}(n\pi r^3)^{\frac13}.
$
So
\begin{eqnarray*}
A'_2
&=&\frac{1}{n\pi r^3} \int_{\Gamma'_2(x)}n\Pr\{Z_3=k-1\}\mathrm{d}y\\
&\leq&\frac{1}{n\pi r^3} \int_{\Gamma'_2(x)}\frac{n(n\pi r^3)^{k-1}e^{-\frac 12\left(n\pi r^3\right)^{\frac13} }}{(k-1)!} \mathrm{d}y\\
&\leq&\frac{1}{n\pi r^3}
\frac{(n\pi r^3)^{k}e^{-\frac {1}{2}\left(n\pi r^3\right)^{\frac13} }}{(k-1)!}
 \sim \frac{o(1)}{n\pi r^3}=o(1).
\end{eqnarray*}
Currently, we have
\begin{eqnarray*}
	 \frac{1}{n\pi r^3} \int_{\Omega\cap B(x,r)}n\Pr\{
	 Z_3=k-1\}\mathrm{d}y\leq A_1+A_2+A_3  =o(1).
\end{eqnarray*}
So Proposition~\ref{proposition1} or (\ref{eq:Pro-5})   holds
under   condition~$C_\mathrm{II}$.
Finally,
\begin{equation}\label{eq:Pro-2-CII-b}
\begin{aligned}
			&\frac{1}{n\pi r^3}  \int_{\Omega\backslash\mathrm{Tube}(\sigma)}n\psi^{k}_{n,r}(x) \mathrm{d}x
		  \int_{\Omega\cap B(x,r)}n\Pr\{Z_3=k-1\}\mathrm{d}y\\
\leq & o(1)  \int_{\Omega}n\psi^{k}_{n,r}(x) \mathrm{d}x\sim o(1)e^{-c}=o(1).
\end{aligned}
\end{equation}

Combining (\ref{eq:Pro-2-CII-a}) and (\ref{eq:Pro-2-CII-b}),
$$
\frac{1}{n\pi r^3}  \int_{\Omega}n\psi^{k}_{n,r}(x) \mathrm{d}x
		  \int_{\Omega\cap B(x,r)}n\Pr\{Z_3=k-1\}\mathrm{d}y=o(1).
$$
This completes the proof of Proposition~\ref{pro:conclusion-1} under condition $C_\mathrm{II}$.

Follow the notes given in Remark~\ref{remark-2:Condition3}, Proposition~\ref{pro:conclusion-1} under condition $C_\mathrm{III}$ can be similarly proved.

\section{Proof of Proposition~\ref{pro:conclusion-2}}\label{sec:Pro-3}

The proof of Proposition~\ref{pro:conclusion-2} is rather simple.  Notice that the assumptions of  Theorem~\ref{thm:Main3D}, i.e., conditions $C_{\mathrm{I}}$ (see Lemma~\ref{lemma:C2-Boundary-volume}), $C_{\mathrm{II}}$ and $C_{\mathrm{III}}$,
allow that
$|B(x,r)\cap\Omega|$ appeared in
$$
\psi^k_{n,r}(x)= \frac{\left(n|B(x,r)\cap\Omega|\right)^k
	e^{-n|B(x,r)\cap\Omega|}}{k!}
$$
to be bounded by $C|B(x,r)|\leq|B(x,r)\cap\Omega|\leq |B(x,r)|$ with some constant $C>0$.

Then follow the same technique demonstrated in~\cite{Penrose-k-connectivity}, or the proof precdure of Conclusion~2 presented in~\cite{Ding2025-IEEETON}, we know that
$$\lim_{n\rightarrow\infty}\Pr\left\{\rho(\chi_n;\delta\geq k+1)=\rho(\chi_n;\kappa\geq
k+1)\right\}=1.
$$
Since Proposition~\ref{pro:conclusion-2} is proved, the proof of Theorem~\ref{thm:Main3D} is therefore completed.

\section{Conclusions and Discussions}
We have demonstrated the two types of critical transmission radii defined in terms of $k-$connectivity and the minimum vertex degree in the random geometry graphs over three-dimensional regions.  The  techniques developed in this article (which are traced back to~\cite{Ding2025-IEEETON} and~\cite{Ding2018-2DRadii}) can also apply to deal with general $d\geq 3$ dimensional    regions.

For any $t\in[0,r]$,   define
\begin{equation*}\label{eq:a(r,t)}
\begin{aligned}
	&a(t)\\
=&|\{(x_1,x_2,\cdots,x_d)\in\mathbb{R}^d: x_1^2+x_2^2+\cdots+x_d^2\leq r^2, x_1\leq t\}|,
\end{aligned}
\end{equation*}
the volume of a major spherical segment in $\mathbb{R}^d$. Let $r=r_n=\left( \frac{\log n+ \frac{dk-d+1}{d-1} \log\log n+\xi}{ \frac{d}{2(d-1)}V_d(1) n }\right)^{\frac1d}$ and $k\geq 0$, then as $n\rightarrow\infty$ or $r\rightarrow 0$,
$$
n\int_0^{\frac r2}\frac{(na(t))^ke^{-na(t)}}{k!}dt\sim
\frac{\left(\frac{d-1}{d}\right)^k\left(\frac{d}{2(d-1)}V_d(1)\right)^{\frac{d-1}{d}}}{e^{\frac{d-1}{d}\xi}{V_{d-1}(1)}k!}.
$$
Here $V_d(1)$ denotes the  volume of   unit $d$-dimensional sphere. Notation $\mathrm{Area}(\partial\Omega)$ denotes the area of hyper surface $\partial\Omega\subset\mathbf{R}^d$.
Based on that, we can similarly prove the following conclusion using the arguments presented in this paper.
\begin{theorem} \label{thm:Main-HD}
	Suppose $\Omega\subset \mathbb{R}^d (d\geq 3) $  is a simply connected compact region with unit-volume.
	The boundary $\partial\Omega$ of the region is $C^2$ smooth and  has finite  area $\mathrm{Area}(\partial\Omega)$, $k\geq 0$ is an integer and $c>0$ is a constant. Point process $\chi_n$ is uniformly distributed on $\Omega$.
	Let
	\begin{equation*}\label{eq:Theorem-radius-HD}
		r_n= \left( \frac{\log n+ \frac{dk-d+1}{d-1} \log\log n+\xi}{ \frac{d}{2(d-1)}V_d(1) n }\right)^{\frac1d},
	\end{equation*}
	where $\xi$ solves
	\begin{equation*}\label{eq:Theorem-radius-1}		 \mathrm{Area}(\partial\Omega)\frac{\left(\frac{d-1}{d}\right)^k\left(\frac{d}{2(d-1)}V_d(1)\right)^{\frac{d-1}{d}}}{e^{\frac{d-1}{d}\xi}{V_{d-1}(1)}k!}  =e^{-c}.
	\end{equation*}
Then  probabilities of the two events
	$\rho(\chi_n;\delta\geq k+1)\leq r_n$ and $\rho(\chi_n;\kappa\geq
	k+1)\leq r_n$ both converge to $\exp\left(-e^{-c}\right)$  as $n\rightarrow\infty$.
\end{theorem}
This conclusion is consistent with Theorem~\ref{thm:Main3D} when $d=3$. We do not discuss the proof of Theorem~\ref{thm:Main-HD} in this paper. Please see manuscript~\cite{Ding2025-HDRadii} for details.


\bibliographystyle{IEEEtran}
\bibliography{3DRGG}

 \end{CJK*}
\end{document}